\documentclass[a4paper,10pt]{amsart}

\usepackage{enumerate, booktabs}
\usepackage{amsmath, amsfonts, amssymb, amsthm}
\usepackage[height=22.5cm, width=16cm]{geometry}

\theoremstyle{plain}
\newtheorem{thm}{Theorem}[section]
\newtheorem{prop}[thm]{Proposition}
\newtheorem{lem}[thm]{Lemma}

\theoremstyle{definition}
\newtheorem{defn}[thm]{Definition}

\theoremstyle{remark}
\newtheorem*{rem}{Remark}

\numberwithin{equation}{section}

\newcommand{\mbb}[1]{\mathbb{#1}}
\newcommand{\lie}[1]{{\mathfrak{#1}}}
\newcommand{\cal}[1]{{\mathcal{#1}}}
\newcommand{\abs}[1]{\lvert #1\rvert}
\newcommand{\norm}[1]{\lVert #1\rVert}

\newcommand{\wt}[1]{\widetilde{#1}}

\DeclareMathOperator{\im}{Im}
\DeclareMathOperator{\re}{Re}

\DeclareMathOperator{\Ad}{Ad}
\DeclareMathOperator{\Int}{Int}

\DeclareMathOperator{\Aut}{Aut}
\DeclareMathOperator{\Lie}{Lie}

\DeclareMathOperator{\Fix}{Fix}
\DeclareMathOperator{\rk}{rk}

\DeclareMathOperator{\codim}{codim}

\DeclareMathOperator{\Span}{Span}

\def\cprime{$'$}
\providecommand{\bysame}{\leavevmode\hbox to3em{\hrulefill}\thinspace}
\providecommand{\MR}{\relax\ifhmode\unskip\space\fi MR }

\providecommand{\href}[2]{#2}

\title{Geometry of invariant domains in complex semi-simple Lie groups}
\author{Christian Miebach}

\begin{document}

\maketitle

\begin{abstract}
We investigate the joint action of two real forms of a semi-simple complex Lie
group $U^\mbb{C}$ by left and right multiplication. After analyzing the orbit
structure, we study the CR structure of closed orbits. The main results are an
explicit formula of the Levi form of closed orbits and the determination of the
Levi cone of generic orbits. Finally, we apply these results to prove
$q$--completeness of certain invariant domains in $U^\mbb{C}$.
\end{abstract}

\section{Introduction}

Let $U^\mbb{C}$ be a connected semi-simple complex Lie group with compact real
form $U$ which is given by the Cartan involution $\theta$. Let us assume that
there are two anti-holomorphic involutive automorphisms $\sigma_1$ and
$\sigma_2$ of $U^\mbb{C}$ which both commute with $\theta$ and let $G_j=
\Fix(\sigma_j)$, $j=1,2$, be the corresponding real forms of $U^\mbb{C}$. The
group $G_1\times G_2$ acts on $U^\mbb{C}$ by $(g_1,g_2)\cdot z:=g_1zg_2^{-1}$.
In this paper we investigate complex-analytic properties of certain $(G_1\times
G_2)$--invariant domains in $U^\mbb{C}$ through the intrinsic Levi form of
closed $(G_1\times G_2)$--orbits.

If $\sigma_1=\sigma_2=\theta$, then we discuss the $(U\times U)$--action on
$U^\mbb{C}$ by left and right multiplication. Every $(U\times U)$--orbit
intersects the set $\exp(i\lie{t})$ in the orbit of the Weyl group
$W:={\cal{N}}_U(\lie{t})/{\cal{Z}}_U(\lie{t})$ where $\lie{t}$ is a
maximal torus in $\lie{u}$. In~\cite{Las1} Lassalle showed that every
bi-invariant domain $\Omega\subset U^\mbb{C}$ is of the form $U\exp(i\omega)U$
for a $W$--invariant domain $\omega\subset\lie{t}$ and that $\Omega$ is a
domain of holomorphy if and only if $\omega$ is convex. In~\cite{AzL} Azad and
Loeb proved the stronger statement that a $(U\times U)$--invariant function
$\Phi$ on $\Omega$ is plurisubharmonic if and only if the $W$--invariant
function
\begin{equation*}
\varphi\colon\lie{t}\to\mbb{R},\quad\varphi(\eta):=\Phi\bigl(\exp(i\eta)\bigr),
\end{equation*}
is convex.

In the case that $\sigma_1=\sigma_2$ and $G_1$ is a real semi-simple Lie group
of Hermitian type there is a distinguished $(G_1\times
G_1)$--invariant in $U^\mbb{C}$, namely the open complex Ol'shanski{\u\i}
semi-group. According to a result of Neeb (\cite{Ne1}) the open Ol'shanski{\u\i}
semi-group is a domain of holomorphy.

Although the above results are statements about complex-analytic properties of
domains in complex Stein manifolds, the method of their proofs is
representation-theoretic. A different approach to the study of
$(G_1\times G_1)$--invariant domains in $U^\mbb{C}$ by analytic methods was made
by Fels an Geatti in~\cite{FeGe}. There, Fels and Geatti gave explicit formulas
for the intrinsic Levi form of a closed orbit $M_z:=(G_1\times G_1)\cdot z$ of
maximal dimension in $U^\mbb{C}$ (in the following called a generic orbit) and
determined the Levi cone of $M_z$, which enabled them to decide whether or not
there may exist a bi-invariant domain of holomorphy containing $z$ in its
boundary.

The main results in this paper are an explicit formula for the intrinsic Levi
form of an arbitrary closed $(G_1\times G_2)$--orbit in $U^\mbb{C}$ and the
determination of the Levi cone of a generic orbit. We use a theorem of Matsuki
(\cite{Ma}) in order to obtain a parameterization of closed $(G_1\times
G_2)$--orbits by certain Cartan algebras in the Lie algebra $\lie{u}^\mbb{C}=
\Lie(U^\mbb{C})$. More precisely, there are finitely many Cartan algebras
$\lie{c}_j$ such that the closed orbits are precisely those intersecting a set
of the form $C_j=n\exp(i\lie{c}_j)$, where the element $n$ can be chosen from
a fixed torus in $U$. It turns out that the weight space decomposition of
$\lie{u}^\mbb{C}$ with respect to $\lie{c}_j$ is well-suited to describe the CR
structure of closed orbits intersecting $C_j$. In particular, the complex
tangent space of such an orbit can be identified with a direct sum of weight
spaces and the intrinsic Levi form of a closed orbit is determined by the Lie
bracket of certain weight vectors together with a coefficient which depends on
the intersection of the orbit with $C_j$ (Theorem~\ref{LeviForm:Formula}). From
this fact it can be derived that the CR structures of closed orbits which belong
to the same set $C_j$ have very similar properties.

The method used here for the derivation of explicit formulas for the Levi form
is different from the one used in~\cite{FeGe}. While Fels and Geatti found
explicit local extensions of complex tangent vectors to CR vector fields on a
generic orbit and computed their Lie brackets, the approach used here avoids
these technical difficulties by pulling back the CR structure of the orbit into
the Lie algebra of $G_1\times G_2$ where the Levi form can be determined via
Lie-theoretic methods. In particular, we obtain a new proof for their results
in the case $G_1=G_2$.

A finer analysis of the weight space decomposition of $\lie{u}^\mbb{C}$ with
respect to $\lie{c}_j$ reveals that it has properties very close to a root
space decomposition. The most important one is the existence of
$\lie{sl}(2)$--triples which enables us to determine the Levi cone of generic
orbits by essentially the same method as in~\cite{FeGe}
(Theorem~\ref{LeviCone}).

In Section~4 we give several applications of the results obtained so far. First
we use the criterion from~\cite{FeGe} together with the knowledge of the Levi
cone in order to decide which $(G_1\times G_2)$--invariant domains containing a
generic orbit in their boundary can be Stein. Secondly, we classify and study
the rank one case in some detail since this case provides a class of examples
where the methods and results become most transparent. This is due to the facts
that complex-analytic properties of smooth domains in Stein manifolds are
determined by the classical Levi form of their boundaries and that in the
rank one case the boundaries of almost all invariant domains coincide with
orbits of hypersurface type. Finally, we define a $(G_1\times G_2)$--invariant
domain $\Omega\subset U^\mbb{C}$ which is the right analogon of the open complex
Ol'shanski{\u\i} semi-group in the case $G_1\not=G_2$. We prove that the
classical Levi form of a $(G_1\times G_2)$--invariant smooth function at a point
$z\in\Omega$ splits into a contribution coming from the complex tangent space of
$(G_1\times G_2)\cdot z$ and a contribution due to a transversal slice. Via this
splitting we construct a strictly $q$--convex function on $\Omega$ which goes to
infinity at $\partial\Omega$, and hence conclude that $\Omega$ is $q$--complete.

\subsection*{Notation}

If $\varphi$ is an automorphism of a Lie group $G$, then by abuse of notation we
write $\varphi$ also for the derived automorphism of $\lie{g}=\Lie(G)$.

\subsection*{Acknowledgment}

This paper is a modified version of the author's Ph.D. thesis~\cite{Mie}. The
support by a Promotionsstipendium of the Studienstiftung des deutschen Volkes
and by SFB/TR~12 of the DFG is gratefully acknowledged.

\section{The $(G_1\times G_2)$--Action on $U^\mbb{C}$}

\subsection{Compatible real forms}

Let $U$ be a connected semi-simple compact Lie group. Then its universal
complexification $U^\mbb{C}$ is a connected semi-simple complex Lie group,
and hence carries a unique structure of a linear
algebraic group (compare~\cite{Chev}). The map $\Phi\colon U\times i\lie{u}\to
U^\mbb{C}$, $(u,\xi)\mapsto u\exp(\xi)$, is a real-analytic diffeomorphism,
called the Cartan decomposition of $U^\mbb{C}$. Furthermore, the map
$\theta\colon U^\mbb{C}\to U^\mbb{C}$,
$\theta\bigl(u\exp(\xi)\bigr):=u\exp(-\xi)$, is an anti-holomorphic involutive
automorphism with $U=\Fix(\theta)$, called the Cartan involution of $U^\mbb{C}$
corresponding to the compact real form $U$. Proofs of these facts can be found
e.\,g.\ in~\cite{Kn}.

Let $\sigma_1$ and $\sigma_2$ be two anti-holomorphic involutive automorphisms
of $U^\mbb{C}$, which both commute with the Cartan involution $\theta$.
The fixed point set $G_j:=\Fix(\sigma_j)$ is a real form of $U^\mbb{C}$ for
$j=1,2$. The assumption that $\sigma_j$ commutes with $\theta$ implies that the
Cartan decomposition of $U^\mbb{C}$ restricts to a real-analytic diffeomorphism
$K_j\times\lie{p}_j\to G_j$, where $K_j:=G_j\cap U$ and $\lie{p}_j
:=\lie{g}_j\cap i\lie{u}$ hold. Thus the real form $G_j$ is a compatible
subgroup of $U^\mbb{C}$ in the sense of~\cite{HeSchw}. In particular,
$K_j$ is a deformation retract of $G_j$.

\begin{rem}
Since $G_j$ is closed, the group $K_j$ is compact and hence a maximal compact
subgroup of $G_j$. Thus $G_j$ has only finitely many connected
components. If the group $U^\mbb{C}$ is simply-connected, it follows
from~\cite{Stei} that $G_j$ is connected.
\end{rem}

The product group $G_1\times G_2$ acts on $U^\mbb{C}$ by left and right
multiplication, i.\,e.\ we define
\begin{equation*}
(g_1,g_2)\cdot z:=g_1zg_2^{-1}
\end{equation*}
where $g_j\in G_j$ and $z\in U^\mbb{C}$.

\begin{defn}\label{DefnOfRegularity}
We say that an element $z\in U^\mbb{C}$ is regular (with respect to
$G_1\times G_2$) if the orbit $(G_1\times G_2)\cdot z$ has maximal dimension.
The element $z$ is called strongly regular (with respect to $G_1\times G_2$), if
it is regular and if $(G_1\times G_2)\cdot z$ is closed. We write $U^\mbb{C}_r$
and $U^\mbb{C}_{sr}$ for the sets of regular and strongly regular elements,
respectively. Finally, we call the orbit $(G_1\times G_2)\cdot z$ generic
if $z$ is strongly regular.
\end{defn}

\begin{rem}
\begin{enumerate}[(a)]
\item If we consider the action of $U^\mbb{C}$ on itself given by conjugation,
then Definition~\ref{DefnOfRegularity} yields the usual notion of (strongly)
regular elements in linear algebraic groups (compare~\cite{Hum3}).
\item The subsets $U^\mbb{C}_r$ and $U^\mbb{C}_{sr}$ are invariant under
$G_1\times G_2$. The set $U^\mbb{C}_{sr}$ can be proven to be open and
dense in $U^\mbb{C}$ which justifies the terminology ``generic orbit''.
\item In~\cite{Ma} an element $z\in U^\mbb{C}$ is called  regular semi-simple
if the automorphism $\Ad(z^{-1})\sigma_1\Ad(z)\sigma_2$ is semi-simple
and if the Lie algebra $\lie{g}_2\cap\Ad(z^{-1})\lie{g}_1$ is Abelian. It can
be shown that an element is regular semi-simple in Matsuki's sense if and only
if it is strongly regular.
\end{enumerate}
\end{rem}

\subsection{The isotropy representation}

The following proposition is crucial. For convenience of the reader we give a
short proof which makes use of the complex-analytic structure of $U^\mbb{C}$.

\begin{prop}\label{CompleteReducibility}
Let $z\in U^\mbb{C}$ be a point such that the orbit $M_z:=(G_1\times G_2)\cdot
z$ is closed. Then the isotropy group $(G_1\times G_2)_z$ is real-reductive
and the isotropy representation of $(G_1\times G_2)_z$ on $T_zU^\mbb{C}$ is
completely reducible.
\end{prop}

\begin{proof}
Since $U^\mbb{C}$ is a Stein manifold, there exists a smooth strictly
plurisubharmonic exhaustion function $\rho\colon U^\mbb{C}\to\mbb{R}$. By
compactness of $U$ we can average $\rho$ using the Haar measure and
hence assume that $\rho$ is $(U\times U)$--invariant. It follows that
$\omega:=i\partial\overline{\partial}\rho$ is a $(U\times U)$--invariant
K\"ahler form on $U^\mbb{C}$ with respect to which $U\times U$ acts in a
Hamiltonian fashion. The last statement means that there exists a $(U\times
U)$--equivariant momentum map $\mu\colon U^\mbb{C}\to\lie{u}^*\oplus\lie{u}^*$.
Since the group $G_1\times G_2$ is compatible with the Cartan decomposition of
$U^\mbb{C}\times U^\mbb{C}$, we can restrict $\mu$ to the subspace
$(i\lie{p}_1)^*\oplus(i\lie{p}_2)^*$ and obtain the restricted momentum map
$\mu_{i\lie{p}}\colon U^\mbb{C}\to(i\lie{p}_1)^*\oplus(i\lie{p}_2)^*$.
According to~\cite{HeSchw} this restricted momentum map encodes a lot of
information about the $(G_1\times G_2)$--action on $U^\mbb{C}$ from which we
need the following.
\begin{enumerate}[(a)]
\item A $(G_1\times G_2)$--orbit is closed in $U^\mbb{C}$ if and only if it
intersects $\cal{M}_{i\lie{p}}:=\mu_{i\lie{p}}^{-1}(0)$ non-trivially
(Proposition~11.2 in~\cite{HeSchw}).
\item If $z\in\cal{M}_{i\lie{p}}$, then the isotropy group $(G_1\times G_2)_z$
is a compatible subgroup of $U^\mbb{C}\times U^\mbb{C}$ and hence
real-reductive (Lemma~5.5 in~\cite{HeSchw}). Together with the previous
statement this implies that isotropy groups of closed orbits are real-reductive.
\item If $z\in\cal{M}_{i\lie{p}}$, then the isotropy representation is
completely reducible (Corollary~14.9 in~\cite{HeSchw}).
\end{enumerate}
Hence, the proposition is proven.
\end{proof}

In the rest of this subsection we will have a closer look at the isotropy
representation. Every element $(\xi_1,\xi_2)\in\lie{g}_1\oplus\lie{g}_2$
induces the tangent vector
\begin{equation}\label{DiffOrbitMap}
\begin{split}
\left.\frac{d}{dt}\right|_0\bigl(\exp(t\xi_1)z\exp(-t\xi_2)\bigr)
&=\left.\frac{d}{dt}\right|_0\bigl(z\exp(t\Ad(z^{-1})\xi_1)\exp(-t\xi_2)\bigr)\\
&=(\ell_z)_*\bigl(\Ad(z^{-1})\xi_1-\xi_2\bigr)\in T_zU^\mbb{C},
\end{split}
\end{equation}
where $\ell_z$ denotes left multiplication with $z\in U^\mbb{C}$. These tangent
vectors span the tangent space of the $(G_1\times G_2)$--orbit through $z$,
i.\,e.\ we obtain $T_zM_z=(\lie{g}_1\oplus\lie{g}_2)\cdot z=\bigl\{(\ell_z)_*
\xi;\ \xi\in\lie{g}_2+\Ad(z^{-1})\lie{g}_1\bigr\}$.

Let $\rho$ denote the isotropy representation of $(G_1\times G_2)_z$
on $T_zU^\mbb{C}$. One checks directly that the isotropy group at $z\in
U^\mbb{C}$ is given by
\begin{equation*}
(G_1\times G_2)_z=\bigl\{(zg_2z^{-1},g_2);\ g_2\in G_2\cap z^{-1}G_1z\bigr\}.
\end{equation*}
Consequently, we may identify $(G_1\times G_2)_z$ with $G_2\cap z^{-1}G_1z$ via
the isomorphism $\Phi\colon G_2\cap z^{-1}G_1z\to(G_1\times G_2)_z$,
$g\mapsto(zgz^{-1},g)$. Similarly, we will identify the tangent space
$T_zM_z$ with $\lie{g}_2+\Ad(z^{-1})\lie{g}_1$
via $(\ell_z)_*$. We conclude from
\begin{align*}
\rho\bigl(\Phi(g)\bigr)(\ell_z)_*\xi&=\left.\frac{d}{dt}\right|_0
(zgz^{-1},g)\cdot\bigl(z\exp(t\xi)\bigr)\\
&=\left.\frac{d}{dt}\right|_0\bigl(zg\exp(t\xi)g^{-1}\bigr)\\
&=\left.\frac{d}{dt}\right|_0z\exp\bigl(t\Ad(g)\xi\bigr)=(\ell_z)_*\Ad(g)\xi
\end{align*}
that the map $(\ell_z)_*$ intertwines the adjoint representation of $G_2\cap
z^{-1}G_1z$ on $\lie{u}^\mbb{C}$ with the isotropy representation of $(G_1\times
G_2)_z$ on $T_zU^\mbb{C}$ modulo $\Phi$. We summarize our considerations in the
following

\begin{prop}\label{IsotrRepn}
Modulo the isomorphism $\Phi$ the isotropy representation of $(G_1\times G_2)_z$
on $T_zU^\mbb{C}$ is equivalent to the adjoint representation of $G_2\cap
z^{-1}G_1z$ on $\lie{u}^\mbb{C}$.
\end{prop}

\subsection{The orbit structure theorem}

We review the main results of~\cite{Ma} in order to describe the orbit
structure of the $(G_1\times G_2)$--action on $U^\mbb{C}$. A proof of
Matsuki's theorem which relies on the momentum map techniques developed
in~\cite{HeSchw} can be found in~\cite{Mie}.

Let $\lie{a}_0$ be a maximal Abelian subspace of $\lie{p}_1\cap\lie{p}_2$ and
let $\lie{t}_0$ be a maximal torus in the centralizer of $\lie{a}_0$ in
$\lie{k}_1\cap \lie{k}_2$. It follows that $\lie{c}_0:=\lie{t}_0\oplus\lie{a}_0$
is a maximally non-compact $\theta$--invariant Cartan subalgebra of
$\lie{g}_1\cap\lie{g}_2$.

\begin{rem}
By maximality of $\lie{a}_0$ the group $A_0^c:=\exp(i\lie{a}_0)$ is a compact
torus in $U$.
\end{rem}

\begin{defn}
A subset of the form $C=n\exp(i\lie{c})\subset U^\mbb{C}$ is called a
standard Cartan subset, if $n\in A_0^c$ and $\lie{c}=\lie{t}\oplus\lie{a}$ is a
$\theta$--stable Cartan subalgebra of $\lie{g}_2\cap\Ad(n^{-1})\lie{g}_1$ such
that $\lie{t}\supset\lie{t}_0$, $\lie{a}\subset\lie{a}_0$ and
$\dim\lie{c}=\dim\lie{c}_0$ hold. The standard Cartan subset
$C_0:=\exp(i\lie{c}_0)$ is called the fundamental Cartan subset.
\end{defn}

We call two standard Cartan subsets equivalent if there is a generic $(G_1
\times G_2)$--orbit which intersects both non-trivially. Let $\{C_j\}_{j\in J}$
be a complete set of representatives for the equivalence classes. For each $j
\in J$ we define the groups
\begin{align*}
{\cal{N}}_{K_1\times K_2}(C_j)&:=\bigl\{(k_1,k_2)\in K_1\times K_2;\ 
k_1C_jk_2^{-1}=C_j\bigr\},\\
{\cal{Z}}_{K_1\times K_2}(C_j)&:=\bigl\{(k_1,k_2)\in K_1\times K_2;\ 
k_1zk_2^{-1}=z\text{ for all $z\in C_j$}\bigr\},
\end{align*}
and $W_{K_1\times K_2}(C_j):={\cal{N}}_{K_1\times K_2}(C_j)/{\cal{Z}}_{K_1\times
K_2}(C_j)$.

\begin{rem}
The group $W_{K_1\times K_2}(C_j)$ is finite for each $j\in J$.
\end{rem}

\begin{thm}[Matsuki]\label{OrbitStructure}
The set $J$ is finite and we have
\begin{equation*}
U^\mbb{C}_{cl}=\bigcup_{j\in J}G_1C_jG_2\quad\text{and}\quad
U^\mbb{C}_{sr}=\dot\bigcup_{j\in J}G_1(C_j\cap U^\mbb{C}_{sr})G_2,
\end{equation*}
where $U^\mbb{C}_{cl}:=\{z\in U^\mbb{C};\ (G_1\times G_2)\cdot z \text{ is
closed}\}$. Moreover, each generic $(G_1\times G_2)$--orbit intersects $C_j$ in
a $W_{K_1\times K_2}(C_j)$--orbit.
\end{thm}

\begin{rem}
If $G_1=G_2$, then let $\lie{c}_0,\dotsc,\lie{c}_k$ be a complete set of
representatives for the equivalence classes of Cartan subalgebras of
$\lie{g}_1$. We can assume without loss of generality that each $\lie{c}_j$ is
$\theta$--stable. Let $\{n_{j,l}\}$ be a set of representatives for the Weyl
group corresponding to $\lie{c}_j$. It can be shown that the sets
$n_{j,l}\exp(i\lie{c}_j)$ exhaust the equivalence classes of standard Cartan
subsets for the $(G_1\times G_1)$--action on $U^\mbb{C}$. Hence, we obtain
Bremigan's theorem (\cite{Bre}).
\end{rem}

\subsection{The weight space decomposition}

Let $C=n\exp(i\lie{c})$ be a standard Cartan subset. In this subsection we
discuss the weight space decomposition
\begin{equation*}
\lie{u}^\mbb{C}=\bigoplus_{\lambda\in\Lambda}\lie{u}^\mbb{C}_\lambda
\end{equation*}
of $\lie{u}^\mbb{C}$ with respect to the Cartan subalgebra
$\lie{c}\subset\lie{g}_2\cap\Ad(n^{-1})\lie{g}_1$. Here, we have written
$\Lambda=\Lambda(\lie{u}^\mbb{C},\lie{c})$ for the set of weights and
$\lie{u}^\mbb{C}_\lambda$ for the weight space corresponding to the weight
$\lambda$. We say that the weight $\lambda$ is real (respectively imaginary) if
$\lambda\not=0$ and $\lambda(\lie{c})\subset\mbb{R}$ (respectively
$\lambda(\lie{c})\subset i\mbb{R}$) holds. A non-zero weight which is neither
real nor imaginary is called complex. We write $\Lambda_r$, $\Lambda_i$ and
$\Lambda_c$ for the sets of real, imaginary and complex weights, and obtain
\begin{equation*}
\Lambda\setminus\{0\}=\Lambda_r\,\dot\cup\,\Lambda_i\,\dot\cup\,\Lambda_c.
\end{equation*}

\begin{rem}
We extend the weight $\lambda$ by $\mbb{C}$--linearity to the complexified
Cartan algebra $\lie{c}^\mbb{C}$. Since $\lambda(\lie{t})\subset i\mbb{R}$
and $\lambda(\lie{a})\subset\mbb{R}$ hold for all $\lambda\in\Lambda$, we
conclude that the weights are real-valued on $i\lie{t}\oplus\lie{a}$.
\end{rem}

Since $n\in A_0^c$, the automorphism
$\tau_n:=\Ad(n^{-1})\sigma_1\Ad(n)\sigma_2\in\Aut(\lie{u}^\mbb{C})$ is
unitary with respect to the Hermitian inner product $\langle\xi_1,\xi_2\rangle
:=-B_{\lie{u}^\mbb{C}}\bigl(\xi_1,\theta(\xi_2)\bigr)$, where
$B_{\lie{u}^\mbb{C}}$ is the Killing form of $\lie{u}^\mbb{C}$.
Consequently, $\tau_n$ is semi-simple with eigenvalues in the unit circle $S^1$.
Since $\tau_n$ leaves $\lie{c}$ pointwise fixed, each weight space
$\lie{u}^\mbb{C}_\lambda$ is invariant under $\tau_n$. Hence,
following~\cite{Ma} we obtain the finer decomposition
\begin{equation}\label{Decomp}
\lie{u}^\mbb{C}=\bigoplus_{(\lambda,a)\in\wt{\Lambda}}
\lie{u}^\mbb{C}_{\lambda,a},
\end{equation}
where $\lie{u}^\mbb{C}_{\lambda,a}:=\bigl\{\xi\in\lie{u}^\mbb{C}_\lambda;\
\tau_n(\xi)=a\xi\bigr\}$ and $\wt{\Lambda}:=\bigl\{(\lambda,a)\in\Lambda\times
S^1;\ \lie{u}^\mbb{C}_{\lambda,a}\not=\{0\}\bigr\}$. The elements of
$\wt{\Lambda}$ are called the extended weights, and~\eqref{Decomp} is called
the extended weight space decomposition.

\begin{rem}
Since $\lie{c}$ is a Cartan subalgebra of $\lie{g}_2\cap\Ad(n^{-1})\lie{g}_1$,
we conclude $\lie{u}^\mbb{C}_{0,1}=\lie{c}^\mbb{C}$.
\end{rem}

We collect some properties of the extended weight space decomposition in the
following

\begin{lem}\label{WeightVectorProperties}
\begin{enumerate}[(1)]
\item The Cartan involution $\theta$ maps $\lie{u}^\mbb{C}_{\lambda,a}$ onto
$\lie{u}^\mbb{C}_{-\lambda,a^{-1}}$. In particular, if $(\lambda,a)$ is an
extended weight, then $(-\lambda,a^{-1})$ is an extended weight, too.
\item We have $B_{\lie{u}^\mbb{C}}\bigl(\lie{u}^\mbb{C}_{\lambda,a},
\lie{u}^\mbb{C}_{\mu,b}\bigr)=0$ unless
$(\lambda,a)=(-\mu,b^{-1})\in\wt{\Lambda}$.
\item Let $\xi_{\lambda,a}\in\lie{u}^\mbb{C}_{\lambda,a}$ with
$\norm{\xi_{\lambda,a}}=1$ be given and let $\eta_{\lambda,a}:= -\bigl[
\xi_{\lambda,a},\theta(\xi_{\lambda,a})\bigr]$. Then we have
$B_{\lie{u}^\mbb{C}}(\eta_{\lambda,a},\eta)=\lambda(\eta)$ for all $\eta\in
\lie{c}$. In particular, $\eta_{\lambda,a}$ does not depend on the
element $a\in S^1$, i.\,e.\ $\eta_{\lambda,a}=\eta_{\lambda,a'}=:\eta_\lambda$
for all $(\lambda,a),(\lambda,a')\in\wt{\Lambda}$.
\item We have $[\xi_{\lambda,a},\xi]=B_{\lie{u}^\mbb{C}}(\xi_{\lambda,a},\xi)
\eta_{\lambda}$ for all $\xi\in\lie{u}^\mbb{C}_{-\lambda,a^{-1}}$.
\end{enumerate}
\end{lem}

\begin{proof}
In order to prove the first claim let
$\eta=\eta_\lie{t}+\eta_\lie{a}\in\lie{t}\oplus\lie{a}=\lie{c}$ and
$\xi\in\lie{u}^\mbb{C}_{\alpha,\lambda}$ be given and consider
\begin{equation*}
\bigl[\eta,\theta(\xi)\bigr]=\theta\bigl[\theta(\eta),\xi\bigr]=
\theta[\eta_\lie{t}-\eta_\lie{a},\xi]=\theta\bigl(\lambda(\eta_\lie{t})
\xi\bigr)-\theta\bigl(\lambda(\eta_\lie{a})\xi\bigr)=-\lambda(\eta)\theta(\xi).
\end{equation*}
Here we used the facts that $\lambda(\lie{t})\subset i\mbb{R}$ while
$\lambda(\lie{a})\subset\mbb{R}$ and that $\theta$ is $\mbb{C}$--anti-linear.
Since $\theta$ commutes with $\tau_n$, we conclude
\begin{equation*}
\tau_n\theta(\xi)=\theta\tau_n(\xi)=\theta(a\xi)=\overline{a}\theta(\xi)=a^{-1}
\theta(\xi),
\end{equation*}
which proves the first claim.

The second claim follows from the fact that the Killing form
$B_{\lie{u}^\mbb{C}}$ is invariant under $\Aut(\lie{u}^\mbb{C})$.

In order to prove the third one we compute
\begin{align*}
B_{\lie{u}^\mbb{C}}(\eta_{\lambda,a},\eta)=-B_{\lie{u}^\mbb{C}}\bigl(
[\xi_{\lambda,a},\theta(\xi_{\lambda,a})],\eta\bigr)&=B_{\lie{u}^\mbb{C}}
\bigl(\theta(\xi_{\lambda,a}),[\xi_{\lambda,a},\eta]\bigr)\\
&=-\lambda(\eta)B_{\lie{u}^\mbb{C}}\bigl(\xi_{\lambda,a},\theta(\xi_{\lambda,a})
\bigr)=\lambda(\eta)\norm{\xi_{\lambda,a}}^2=\lambda(\eta).
\end{align*}

The last claim is proven in the same way as Lemma~2.18(a) in~\cite{Kn}.
\end{proof}

Standard arguments from Lie theory (see for example Chapter~II.4
in~\cite{Kn}) lead to the following result.

\begin{prop}\label{WeightSpaceDecomposition}
\begin{enumerate}[(1)]
\item Let $\lambda\not=0$. After a suitable normalization the elements
$\eta_\lambda$, $\xi_{\lambda,a}$ and $\theta(\xi_{\lambda,a})$ form an
$\lie{sl}(2)$--triple.
\item If $\lambda\not=0$, then we have $\dim_\mbb{C}\lie{u}^\mbb{C}_{\lambda,a}
=1$ and $\dim_\mbb{C}\lie{u}^\mbb{C}_{m\lambda,a^m}=0$ for all $m\geq2$.
\item The set $\Lambda\setminus\{0\}$ of non-zero weights fulfills the axioms
of an abstract root system in $(i\lie{t}\oplus\lie{a})^*$.
\item Let $\lambda,\mu\in\Lambda\setminus\{0\}$ such that
$\lambda+\mu\in\Lambda\setminus\{0\}$ holds. Then we have
$[\lie{u}^\mbb{C}_\lambda,\lie{u}^\mbb{C}_\mu]=\lie{u}^\mbb{C}_{\lambda+\mu}$.
\end{enumerate}
\end{prop}

\section{CR Geometry of Closed Orbits}

\subsection{Preliminaries from CR geometry}

In this subsection we will review the basic definitions and facts from the
theory of CR submanifolds as far as they are needed later on. For more details
and complete proofs we refer the reader to the textbooks~\cite{BER}
and~\cite{Bog}.

Let $Z$ be a complex manifold with complex structure $J$. A real submanifold $M$
of $Z$ is called a Cauchy-Riemann or CR submanifold if the dimension of the
complex tangent space $H_pM:=T_pM\cap J_pT_pM$ does not depend on the point
$p\in M$. In this case, the set $HM:=\bigcup_{p\in M}H_pM$ is a smooth subbundle
of the tangent bundle $TM$ invariant under the complex structure $J$, called
the complex tangent bundle of $M$. A CR submanifold $M\subset Z$ is called
generic if $T_pM+J_pT_pM=T_pZ$ holds for all $p\in M$. For example, every
smooth real hypersurface in $Z$ is a generic CR submanifold of $Z$.

\begin{rem}
Since the group $G_1\times G_2$ acts by holomorphic transformations on
$U^\mbb{C}$, each closed $(G_1\times G_2)$--orbit is a CR submanifold of
$U^\mbb{C}$. Since the $(G_1\times G_2)$--action extends to a transitive
$(U^\mbb{C}\times U^\mbb{C})$--action on $U^\mbb{C}$, each closed orbit is
moreover generic as a CR submanifold.
\end{rem}

A smooth section in $HM$ is called a CR vector field on $M$. A smooth map $f$
from $M$ into a CR submanifold $M'\subset(Z',J')$ is called a CR map if $f_*$
maps $HM$ into $HM'$ and if $f_*J=J'f_*$ holds. A CR function on $M$ is a CR map
$M\to\mbb{C}$, where $\mbb{C}$ is equipped with its usual structure as complex
manifold.

For each CR submanifold $M\subset Z$ one can define the intrinsic Levi form,
which generalizes the classical Levi form of a smooth hypersurface.

\begin{defn}
The Levi form of $M$ at the point $p$ is the map $\cal{L}_p\colon H_pM\times
H_pM\to T^\mbb{C}_pM/H^\mbb{C}_pM$ defined by
\begin{equation*}
\cal{L}_p(v,w):=\left(\frac{i}{2}[V,W]_p-\frac{1}{2}[V,JW]_p\right)\bmod
H^\mbb{C}_pM,
\end{equation*}
where $V$ and $W$ are CR vector fields on $M$ with $V_p=v$ and $W_p=w$.
\end{defn}

\begin{rem}
One can show that the intrinsic Levi form is well-defined, i.\,e.\ that it does
not depend on the choice of CR extensions of $v,w\in H_pM$ (compare~\cite{Bog}).
\end{rem}

The Levi cone $\cal{C}_p$ of $M$ at $p$ is by definition the closed convex cone
generated by the vectors $\cal{L}_z(v,v)$ where $v$ runs through $H_pM$. Because
of $\cal{L}_p(v,w)=\overline{\cal{L}_p(w,v)}$ the Levi cone is contained in
$T_pM/H_pM$. The Levi cone generalizes the signature of the classical Levi form
of a hypersurface. Its significance stems from the fact that it governs the
local extension of CR functions on $M$ to holomorphic functions on $Z$.

\begin{thm}[Boggess, Polking]\label{BoggessPolking}
Let $M$ be a generic CR submanifold of a complex manifold $Z$ and let us assume
that the Levi cone at some point $p\in M$ satisfies
${\cal{C}}_p(M)=T_pM/H_pM$. Then, for each neighborhood $\omega$ of $p$ in $M$
there exists a neighborhood $\Omega$ of $p$ in $Z$ satisfying $\Omega\cap
M\subset\omega$ which has the property that every CR function on $\Omega\cap M$
extends to a unique holomorphic function on $\Omega$.
\end{thm}

A proof of this theorem can be found in~\cite{Bog}.

\subsection{The complex tangent space of a closed orbit}

Let $z\in U^\mbb{C}$ be given such that the orbit $M_z=(G_1\times G_2)\cdot z$
is closed in $U^\mbb{C}$. By Matsuki's theorem we can assume that there is a
standard Cartan subset $C=n\exp(i\lie{c})$ which contains $z=n\exp(i\eta)$.
We define
\begin{equation*}
\wt{\Lambda}(z):=\bigl\{(\lambda,a)\in\wt{\Lambda};\ ae^{-2i\lambda(\eta)}
=1\bigr\}
\end{equation*}
and set $\tau_z:=\Ad(z^{-1})\sigma_1\Ad(z)\sigma_2\in\Aut(\lie{u}^\mbb{C})$.

\begin{lem}\label{tau_z}
The automorphism $\tau_z$ is semi-simple and we have
\begin{equation*}
\Fix(\tau_z)=\bigl(\lie{g}_2\cap\Ad(z^{-1})\lie{g}_1\bigr)^\mbb{C}
=\bigoplus_{(\lambda,a)\in\wt{\Lambda}(z)}\lie{u}^\mbb{C}_{\lambda,a}.
\end{equation*}
\end{lem}

\begin{proof}
The first equality is a consequence of~\cite{Ma}, p.~57. In order to prove the
second one let $\xi=\sum_{(\lambda,a)}\xi_{\lambda,a}$ be an arbitrary element
of
$\lie{u}^\mbb{C}$. Then we have
\begin{align*}
\tau_z(\xi)=\Ad(z^{-1})\sigma_1\Ad(z)\sigma_2(\xi)&=
\Ad\bigl(\exp(-i\eta)\bigr)\tau_n\Ad\bigl(\exp(-i\eta)\bigr)\xi\\
&=\Ad\bigl(\exp(-i\eta)\bigr)\tau_n\left(\sum_{(\lambda,a)}
e^{-i\lambda(\eta)}\xi_{\lambda,a}\right)\\
&=\Ad\bigl(\exp(-i\eta)\bigr)\sum_{(\lambda,a)}a
e^{-i\lambda(\eta)}\xi_{\lambda,a }=\sum_{(\lambda,a)}a
e^{-2i\lambda(\eta)}\xi_{\lambda,a}.
\end{align*}
This proves that $\tau_z$ is semi-simple. Moreover, $\tau_z(\xi)=\xi$ holds if
and only if $\xi_{\lambda,a}=0$ for all $(\lambda,a)\notin\wt{\Lambda}(z)$.
\end{proof}

Since $\lie{g}_2\cap\Ad(z^{-1})\lie{g}_1$ is isomorphic to the Lie algebra of
$(G_1\times G_2)_z$, we obtain the following characterization of strongly
regular elements in terms of the extended weights as a corollary.

\begin{thm}
We have $\codim_\mbb{R}(G_1\times G_2)\cdot z=
\dim_\mbb{R}\lie{c}+(\#\wt{\Lambda}(z)-1)$. The element $z$ is strongly regular
if and only if $\wt{\Lambda}(z)=\bigl\{(0,1)\bigr\}$ holds. This implies that
the codimension of a generic orbit coincides with the rank of the real-reductive
Lie algebra $\lie{g}_1\cap \lie{g}_2$.
\end{thm}

Finally we describe the tangent space $T_zM_z$ in terms of the extended weight
space decomposition.

\begin{thm}\label{ComplexTangentSpace}
Under the map $(\ell_z)_*$ the tangent space $T_zM_z$ is isomorphic to
\begin{equation*}
\lie{g}_2+\Ad(z^{-1})\lie{g}_1=
\bigl(\lie{g}_2\cap\Ad(z^{-1})\lie{g}_1\bigr)\oplus\bigoplus_{(\lambda,a)\notin
\wt{\Lambda}(z)}\lie{u}^\mbb{C}_{\lambda,a}.
\end{equation*}
In particular, the complex tangent space of $(G_1\times G_2)\cdot z$ is
isomorphic to $\bigoplus_{(\lambda,a)\notin\wt{\Lambda}(z)}
\lie{u}^\mbb{C}_{\lambda,a}$.
\end{thm}

\begin{rem}
From now on we will identify the quotient $T^\mbb{C}_zM/H^\mbb{C}_zM$ with
$R^\mbb{C}_zM:=(\ell_z)_*\bigl(\lie{g}_2\cap\Ad(z^{-1})\lie{g}_1\bigr)^\mbb{C}$.
It follows that these spaces are isomorphic as $(G_1\times G_2)_z$--modules.
\end{rem}

\begin{proof}[Proof of Theorem~\ref{ComplexTangentSpace}]
Since $\tau_z$ is semi-simple, we conclude from Lemma~1(i) in~\cite{Ma} that
\begin{equation*}
\lie{u}^\mbb{C}=i\bigl(\lie{g}_2\cap\Ad(z^{-1})\lie{g}_1\bigr)\oplus\bigl(
\lie{g}_2+\Ad(z^{-1})\lie{g}_1\bigr)
\end{equation*}
holds. Moreover, one checks directly that this decomposition is orthogonal with
respect to the real part of the Killing form $B_{\lie{u}^\mbb{C}}$. Similarly,
we have the decomposition
\begin{equation*}
\lie{u}^\mbb{C}=\Fix(\tau_z)\oplus\Fix(\tau_z)^\perp,\quad
\Fix(\tau_z)=\bigl(\lie{g}_2\cap\Ad(z^{-1})\lie{g}_1\bigr)^\mbb{C},
\end{equation*}
where the orthogonal complement $\Fix(\tau_z)^\perp$ with respect to
$B_{\lie{u}^\mbb{C}}$ is the sum of the $\tau_z$--eigenspaces corresponding to
eigenvalues $\not=1$. These observations imply
\begin{equation*}
\lie{g}_2+\Ad(z^{-1})\lie{g}_1=\bigl(\lie{g}_2\cap\Ad(z^{-1})\lie{g}_1\bigr)
\oplus\Fix(\tau_z)^\perp.
\end{equation*}
Since the same argument as the one in the proof of Lemma~\ref{tau_z} implies the
equality
\begin{equation*}
\Fix(\tau_z)^\perp=\bigoplus_{(\lambda,a)\notin\wt{\Lambda}(z)}
\lie{u}^\mbb{C}_{\lambda,a},
\end{equation*}
the theorem is proven.
\end{proof}

\subsection{Pulling back the Levi form into the Lie algebra}

As abbreviation we put $G:=G_1\times G_2$ in this subsection. Consequently, we
have $\lie{g}:=\Lie(G)=\lie{g}_1\oplus\lie{g}_2$.

As we have remarked above, every closed $G$--orbit $M_z=G\cdot z$ is a generic
CR submanifold of $U^\mbb{C}$. Let $\pi_z\colon\lie{g}\to\lie{g}\cdot z=
T_zM_z$ be the differential of the orbit map. By Equation~\eqref{DiffOrbitMap}
the map $\pi_z$ is given by
\begin{equation*}
\pi_z(\xi_1,\xi_2)=(\ell_z)_*\bigl(\Ad(z^{-1})\xi_1-\xi_2\bigr).
\end{equation*}
In this subsection we will pull back the CR structure of $M_z$ into the Lie
algebra $\lie{g}$ and compute the Levi form of $M_z$ via this pull back. The
following proposition is essential.

\begin{prop}
We have the $G_z$--invariant decomposition $\lie{g}=\lie{g}_z\oplus
\lie{q}_z$, and $\lie{q}_z$ and $T_zM_z$ are isomorphic as $G_z$--spaces where
the isomorphism is given by $\wt{\pi}_z:=\pi_z|_{\lie{q}_z}$. Since the
complex tangent space $H_zM_z$ is invariant under $G_z$, we obtain the
$G_z$--invariant decomposition $\lie{q}_z=R(\lie{q}_z)\oplus H(\lie{q}_z)$ where
$H(\lie{q}_z):=\wt{\pi}_z^{-1}(H_zM_z)$ and $R(\lie{q}_z):=\wt{\pi}_z^{-1}
(R_zM_z)$.
\end{prop}

\begin{proof}
We only have to show that the adjoint representation of $G_z$ on $\lie{g}$ is
completely reducible. This follows from Proposition~\ref{CompleteReducibility}
since $G_z$ is conjugate to a compatible subgroup of $U^\mbb{C}\times
U^\mbb{C}$ if the orbit $G\cdot z=M_z$ is closed.
\end{proof}

\begin{prop}\label{LeviForm}
The Levi form $\cal{L}_z\colon H_zM_z\times H_zM_z\to R^\mbb{C}_zM_z$ is given
by
\begin{equation}
\cal{L}_z(v,w)=\pi_z\left(\frac{i}{2}\left[\wt{\pi}_z^{-1}(v),\wt{\pi}_z^{-1}
(w)\right]- \frac{1}{2}\left[\wt{\pi}_z^{-1}(v),\wt{\pi}_z^{-1}(iw)\right]
\right)\bmod H^\mbb{C}_zM_z.
\end{equation}
\end{prop}

\begin{proof}
Let $v,w\in H_zM_z$ be given and let $V,W$ be CR vector fields on $M_z$ with
$V_z=v$ and $W_z=w$. Since the orbit map $G\to M_z$ is a $G_z$--principal
bundle, there exist projectable vector fields $\wt{V}$ and $\wt{W}$ on $G$
with $\wt{V}_e=\wt{\pi}_z^{-1}(v)$ and $\wt{W}_e=\wt{\pi}_z^{-1}(w)$ such that
$\pi_z\wt{V}=V$ and $\pi_z\wt{W}=W$ hold. For a proof of this fact and more
details about projectable vector fields we refer the reader to~\cite{KobNo1}.
Although it is in general not possible to choose the vector fields $\wt{V}$ and
$\wt{W}$ to be left-invariant, the same argument which proves well-definedness
of the intrinsic Levi form applies to show that
\begin{equation*}
\left(\frac{i}{2}[\wt{V},\wt{W}]_e-\frac{1}{2}[\wt{V},\wt{JW}]_e\right)\bmod
H^\mbb{C}(\lie{q}_z)
\end{equation*}
does only depend on the values $\wt{V}_e$ and $\wt{W}_e$ (compare the proof of
Lemma~1 in Chapter~10.1 of~\cite{Bog}). Therefore we conclude
\begin{equation*}
\left(\frac{i}{2}[\wt{V},\wt{W}]_e-\frac{1}{2}[\wt{V},\wt{JW}]_e\right)\bmod
H^\mbb{C}(\lie{q}_z)=\left(\frac{i}{2}\left[\wt{\pi}_z^{-1}(v),\wt{\pi}_z^{-1}
(w)\right]- \frac{1}{2}\left[\wt{\pi}_z^{-1}(v),\wt{\pi}_z^{-1}(iw)\right]
\right)\bmod H^\mbb{C}(\lie{q}_z),
\end{equation*}
and obtain
\begin{align*}
\cal{L}_z(v,w)&=\left(\frac{i}{2}[V,W]_z-\frac{1}{2}[V,JW]_z\right) \bmod
H^\mbb{C}_zM_z\\
&=\left(\frac{i}{2}[\pi_z\wt{V},\pi_z\wt{W}]_z-\frac{1}{2}[\pi_z\wt{V},
\pi_z\wt{JW}] _z\right) \bmod H^\mbb{C}_zM_z\\
&=\left(\frac{i}{2}\pi_z[\wt{V},\wt{W}]_e-\frac{1}{2}\pi_z[\wt{V},\wt{JW}]_e
\right) \bmod H^\mbb{C}_zM_z\\
&=\pi_z\left(\frac{i}{2}[\wt{V},\wt{W}]_e-\frac{1}{2}[\wt{V},\wt{JW}]_e
\bmod H^\mbb{C}(\lie{q}_z)\right)\\
&=\pi_z\left(\frac{i}{2}\left[\wt{\pi}_z^{-1}(v),\wt{\pi}_z^{-1}
(w)\right]- \frac{1}{2}\left[\wt{\pi}_z^{-1}(v),\wt{\pi}_z^{-1}(iw)\right]
\bmod H^\mbb{C}(\lie{q}_z)\right)\\
&=\pi_z\left(\frac{i}{2}\left[\wt{\pi}_z^{-1}(v),\wt{\pi}_z^{-1}
(w)\right]- \frac{1}{2}\left[\wt{\pi}_z^{-1}(v),\wt{\pi}_z^{-1}(iw)\right]
\right)\bmod H^\mbb{C}_zM_z.
\end{align*}
This finishes the proof.
\end{proof}

In the next subsection we will use the weight space decomposition in order to
determine the map $\wt{\pi}_z^{-1}$ explicitely.

\subsection{The Levi form of a closed orbit}

This rather technical subsection contains the computations which are necessary
to achieve the final formula of the Levi form.

\begin{lem}
We have
\begin{equation*}
\lie{g}_z=\ker(\pi_z)=\bigl\{(\Ad(z)\xi,\xi);\
\xi\in\lie{g}_2\cap\Ad(z^{-1})\lie {g}_1\bigr\}.
\end{equation*}
\end{lem}

\begin{proof}
One checks directly that $\bigl\{(\Ad(z)\xi,\xi);\ \xi\in\lie{g}_2\cap
\Ad(z^{-1})\lie {g}_1\bigr\}\subset\ker(\pi_z)$ holds. The other inclusion
follows for dimensional reasons.
\end{proof}

\begin{lem}\label{InvImage}
The subspace $\lie{q}_z=R(\lie{q}_z)\oplus H(\lie{q}_z)$ is determined by the
following.
\begin{enumerate}[(i)]
\item We have
\begin{equation*}
R(\lie{q}_z)=\wt{\pi}_z^{-1}\bigl(\lie{g}_2\cap\Ad(z^{-1})\lie{g}_1\bigr)
=\bigl\{(\Ad(z)\xi,-\xi);\ \xi\in\lie{g}_2\cap\Ad(z^{-1})\lie{g}_1\bigr\}.
\end{equation*}
\item We have
\begin{equation*}
(\wt{\pi}_z)^{-1}(\lie{u}^\mbb{C}_{\lambda,a})=\bigl\{
\bigl(\Ad(z)\sigma_2(\xi)+\sigma_1\bigl(\Ad(z)\sigma_2(\xi)\bigr),
\xi+\sigma_2(\xi)\bigr);\ \xi\in\lie{u}^\mbb{C}_{\lambda,a}\bigr\}
\end{equation*}
for all $(\lambda,a)\in\wt{\Lambda}\setminus\wt{\Lambda}(z)$.
\end{enumerate}
\end{lem}

\begin{proof}
Firstly, we have to show that $\bigl\{(\Ad(z)\xi,-\xi);\
\xi\in\lie{g}_2\cap\Ad(z^{-1})\lie{g}_1\bigr\}$ is contained in $\lie{q}_z=
\lie{g}_z^\perp$ where the orthogonal complement is taken with respect to the
Killing form of $\lie{g}_1\oplus\lie{g}_2$. Hence, let
$\xi,\xi'\in\lie{g}_2\cap\Ad(z^{-1})\lie{g}_1$ and consider
\begin{equation*}
\begin{split}
B_{\lie{g}_1\oplus\lie{g}_2}\bigl((\Ad(z)\xi,\xi),(\Ad(z)\xi',-\xi')\bigr)
&=B_{\lie{g}_1}(\Ad(z)\xi,\Ad(z)\xi')-B_{\lie{g}_2}(\xi,\xi')\\
&=B_{\lie{u}^\mbb{C}}(\xi,\xi')-B_{\lie{u}^\mbb{C}}(\xi,\xi')=0.
\end{split}
\end{equation*}
A simple computation shows $\bigl\{(\Ad(z)\xi,-\xi);\
\xi\in\lie{g}_2\cap\Ad(z^{-1})\lie{g}_1\bigr\}\subset R(\lie{q}_z)$. Since the
converse inclusion follows for dimensional reasons, the first claim is proven.

A similar argument as above implies that $\bigl\{
\bigl(\Ad(z)\sigma_2(\xi)+\sigma_1\bigl(\Ad(z)\sigma_2(\xi)\bigr),
\xi+\sigma_2(\xi)\bigr);\ \xi\in\lie{u}^\mbb{C}_{\lambda,a}\bigr\}$ lies in
$\lie{q}_z$. In order to prove the second assertion let
$\xi\in\lie{u}^\mbb{C}_{\lambda,a}$ be given and consider
\begin{align*}
\wt{\pi}_z\bigl(\Ad(z)\sigma_2(\xi)+\sigma_1\bigl(\Ad(z)\sigma_2(\xi)\bigr),
\xi+\sigma_2(\xi)\bigr)&=\sigma_2(\xi)+\Ad(z^{-1})\sigma_1\Ad(z)\sigma_2(\xi)
-\xi-\sigma_2(\xi)\\
&=\tau_z(\xi)-\xi=\bigl(ae^{-2i\lambda(\eta)}-1\bigr)\xi=:
\varphi_{\lambda,a}(\xi).
\end{align*}
Since $\varphi_{\lambda,a}(\xi)\in\lie{u}^\mbb{C}_{\lambda,a}$ holds, the lemma
is proven.
\end{proof}

\begin{rem}
Note that the map
$\varphi_{\lambda,a}\colon\lie{u}^\mbb{C}_{\lambda,a}
\to\lie{u}^\mbb{C}_{\lambda,a}$ is an isomorphism if and only if
$(\lambda,a)\notin\wt{\Delta}(z)$ holds. In this case the inverse
map is given by
\begin{equation*}
\varphi_{\lambda,a}^{-1}(\xi)=\frac{1}{ae^{-2i\lambda(\eta)}-1}\xi.
\end{equation*}
\end{rem}

\begin{defn}
A Levi basis of $H_zM_z$ is a basis $(\xi_{\lambda,a})_{(\lambda,a)}$ of
$H_zM_z$ such that $\xi_{\lambda,a}\in\lie{u}^\mbb{C}_{\lambda,a}$
and $\sigma_2(\xi_{\lambda,a})=\xi_{\sigma_2(\lambda),a}$ hold
for all $(\lambda,a)\in\wt{\Lambda}\setminus\wt{\Lambda}(z)$.
\end{defn}

From now on we fix a Levi basis $(\xi_{\lambda,a})$ of $H_zM_z$.

\begin{thm}\label{LeviForm:Formula}
We obtain the following formula for the Levi form of $M_z$:
\begin{equation*}
{\cal{L}}_z(\xi_{\lambda,a},\xi_{\mu,b})=
\begin{cases}
\frac{i}{ae^{-2i\lambda(\eta)}-1}\bigl[\xi_{\lambda,a},
\xi_{\sigma_2(\mu),b}\bigr]&\text{if }
(\lambda+\sigma_2(\mu),ab)\in\wt{\Lambda}(z)\\
0 & \text{else}
\end{cases}.
\end{equation*}
\end{thm}

\begin{proof}
We will start by computing $\wt{\pi}_z^{-1}(\xi_{\lambda,a})$ for
$\xi_{\lambda,a}\in\lie{u}^\mbb{C}_{\lambda,a}$. Lemma~\ref{InvImage}
gives us
\begin{align*}
\wt{\pi}_z^{-1}(\xi)&=\left(\Ad(z)\sigma_2(\varphi_{\lambda,a}^{-1}\xi)
+\sigma_1\bigl(\Ad(z)\sigma_2(\varphi_{\lambda,a}^{-1}\xi)\bigr),
\varphi_{\lambda,a}^{-1}\xi+\sigma_2(\varphi_{\lambda,a}^{-1}
\xi)\right)\\
&=\left(\Ad(z)\bigl(\sigma_2(\varphi_{\lambda,a}^{-1}\xi)
+\tau_z(\varphi_{\lambda,a}^{-1}\xi)\bigr),
\varphi_{\lambda,a}^{-1}\xi+\sigma_2(\varphi_{\lambda,a}^{-1}
\xi)\right)\\
&=\left(\Ad(z)\bigl(\varphi_{\lambda,a}^{-1}\xi
+\sigma_2(\varphi_{\lambda,a}^{-1}\xi)\bigr)+\Ad(z)\xi,
\varphi_{\lambda,a}^{-1}\xi+\sigma_2(\varphi_{\lambda,a}^{-1}
\xi)\right)
\end{align*}
for any $\xi\in\lie{u}^\mbb{C}_{\lambda,a}$. In the next step we determine
the Lie bracket
$\bigl[\wt{\pi}^{-1}_z(\xi_{\lambda,a}),\wt{\pi}^{-1}_z
(\xi_{\mu,b})\bigr]$. Since the Lie bracket of
$\lie{g}=\lie{g}_1\oplus\lie{g}_2$ is defined component-wise, we consider
\begin{equation*}
\begin{split}
\left[\Ad(z)\bigl(\varphi_{\lambda,a}^{-1}\xi_{\lambda,a}+\sigma_2
(\varphi_{\lambda,a}^{-1}\xi_{\lambda,a})\bigr)+\Ad(z)\xi_{\alpha,
\lambda},\
\Ad(z)\bigl(\varphi_{\mu,b}^{-1}\xi_{\mu,b}+\sigma_2
(\varphi_{\mu,b}^{-1}\xi_{\mu,b})\bigr)+\Ad(z)\xi_{\mu,b}\right]\\
=\Ad(z)\left[\varphi_{\lambda,a}^{-1}\xi_{\lambda,a}+\sigma_2
(\varphi_{\lambda,a}^{-1}\xi_{\lambda,a})+\xi_{\lambda,a},\
\varphi_{\mu,b}^{-1}\xi_{\mu,b}+\sigma_2
(\varphi_{\mu,b}^{-1}\xi_{\mu,b})+\xi_{\mu,b}\right]
\end{split}
\end{equation*}
and
\begin{equation*}
\left[\varphi_{\lambda,a}^{-1}\xi_{\lambda,a}+\sigma_2(\varphi_{\alpha,
\lambda}^{-1}\xi_{\lambda,a}),\
\varphi_{\mu,b}^{-1}\xi_{\mu,b}+\sigma_2(\varphi_{\mu,b}^{-1}\xi_{
\mu,b})\right].
\end{equation*}
The application of $\pi_z$ to the element in $\lie{g}_1\oplus\lie{g}_2$ whose
components are given by the above gives
\begin{equation}\label{Summand1}
\left[\varphi_{\lambda,a}^{-1}\xi_{\lambda,a}+\sigma_2(\varphi_{\alpha
,\lambda}^{-1}\xi_{\lambda,a}),
\xi_{\mu,b}\right]+\left[\xi_{\lambda,a},\varphi_{\mu,b}^{-1}\xi_{
\mu,b}+\sigma_2(\varphi_{\mu,b}^{-1}
\xi_{\mu,b})\right]+[\xi_{\lambda,a},\xi_{\mu,b}].
\end{equation}
By the same computation we obtain for
$\pi_z\bigl([\wt{\pi}^{-1}_z(\xi_1),\wt{\pi}^{-1}_z(i\xi_2)]
\bigr)$
the following expression:
\begin{equation}\label{Summand2}
\left[\varphi_{\lambda,a}^{-1}\xi_{\lambda,a}+\sigma_2(\varphi_{\alpha
,\lambda}^{-1}\xi_{\lambda,a}),
i\xi_{\mu,b}\right]+\left[\xi_{\lambda,a},\varphi_{\mu,b}^{-1}i\xi_
{\mu,b}+\sigma_2(\varphi_{\mu,b}^{-1}
i\xi_{\mu,b})\right]+[\xi_{\lambda,a},i\xi_{\mu,b}].
\end{equation}
To arrive at the Levi form, we have to multiply~\eqref{Summand1} by
$\frac{i}{2}$ and subtract~\eqref{Summand2} multiplied
by $\frac{1}{2}$. Due to the facts that $\varphi_{\lambda,a}$ and
$\varphi_{\mu,b}$ are complex-linear, while $\sigma_2$ is anti-linear over
$\mbb{C}$, this leads to
\begin{equation*}
i\bigl[\xi_{\lambda,a},\sigma_2(\varphi_{\mu,b}^{-1}\xi_{\mu,b}
)\bigr].
\end{equation*}
Inserting the concrete expression for $\varphi_{\mu,b}^{-1}$ yields
\begin{equation*}
\pi_z\left(\frac{i}{2}\bigl[\wt{\pi}_z^{-1}(\xi_{\lambda,a}),
\wt{\pi}_z^{-1}(\xi_{\mu,b})\bigr]-\frac{1}{2}\bigl[\wt{\pi}
_z^{-1}
(\xi_{\lambda,a}),\wt{\pi}_z^{-1}(i\xi_{\mu,b})\bigr]\right)=
\frac{i}{\mu^{-1}e^{2i\overline{\beta(\eta)}}-1}\bigl[\xi_{\lambda,a},
\sigma_2(\xi_{\mu,b})\bigr].
\end{equation*}
To arrive at the Levi form, we have to project this element onto
$\bigl(\lie{g}_2\cap\Ad(z^{-1})\lie{g}_1\bigr)^\mbb{C}$. Consequently, we
only obtain a nonzero contribution if
$\bigl[\xi_{\lambda,a},\sigma_2(\xi_{\mu,b})\bigr]
\in\Fix(\tau_z)$ holds. By the definition of a Levi basis this condition
translates into the one formulated in the theorem. This finishes the proof.
\end{proof}

\subsection{The quadratic Levi form}

In this subsection we will derive explicit formulas for the quadratic Levi form
of a generic orbit $M_z=(G_1\times G_2)\cdot z$ from
Theorem~\ref{LeviForm:Formula}. For $(\lambda,a)\in\wt{\Lambda}$ we define
\begin{equation*}
\lie{u}^\mbb{C}[\lambda,a]:=\lie{u}^\mbb{C}_{\lambda,a}+
\lie{u}^\mbb{C}_{\sigma_2(\lambda),a}+\lie{u}^\mbb{C}_{-\lambda,a^{-1 }}
+\lie{u}^\mbb{C}_{-\sigma_2(\lambda),a^{-1}}.
\end{equation*}
Since the $\lie{u}^\mbb{C}[\lambda,a]\perp\lie{u}^\mbb{C}[\mu,b]$  with respect
to the Levi form ${\cal{L}}_z$ for
$\bigl(\lambda+\sigma_2(\mu),ab\bigr)\notin\wt{\Lambda}(z)=\bigl\{(0,1)\bigr\}$,
the Levi form is determined by its restriction to these spaces for which
explicit formulas are given in the next proposition. We make use of the
partition
$\Lambda\setminus\{0\}=\Lambda_r\,\dot\cup\,\Lambda_i\,\dot\cup\,\Lambda_c$.

\begin{prop}\label{QuadrLeviForm}
Let $(\lambda,a)\in\wt{\Lambda}$ be given.
\begin{enumerate}[(1)]
\item For $\lambda\in\Lambda_r$ we obtain
\begin{equation*}
\widehat{\cal{L}}_z(r_\lambda\xi_{\lambda,a}+r_{-\lambda}\xi_{-\lambda,a^{-1}})=
-2\im\left(\frac{r_\lambda\overline{r}_{-\lambda}}{ae^{-2i\lambda(\eta)}-1}
\right)[\xi_{\lambda,a},\xi_{-\lambda,a^{-1}}].
\end{equation*}
\item For $\lambda\in\Lambda_i$ and $a=1$ we obtain
\begin{equation*}
\widehat{\cal{L}}_z(r_\lambda\xi_{\lambda,1}+r_{-\lambda}\xi_{-\lambda,1})=
\left(\frac{\abs{r_\lambda}^2}{e^{-2i\lambda(\eta)}-1}-\frac{\abs{r_{-\lambda}}
^2}{e^{2i\lambda(\eta)}-1}\right)i[\xi_{\lambda,1},\xi_{-\lambda,1}].
\end{equation*}
\item For $\lambda\in\Lambda_i$ and $a=-1$ we obtain
\begin{equation*}
\widehat{\cal{L}}_z(r_\lambda\xi_{\lambda,-1}+r_{-\lambda}\xi_{-\lambda,-1})=
-\left(\frac{\abs{r_\lambda}^2}{e^{-2i\lambda(\eta)}+1}-\frac{\abs{r_{-\lambda}}
^2}{e^{2i\lambda(\eta)}+1}\right)i[\xi_{\lambda,-1},\xi_{-\lambda,-1}].
\end{equation*}
\item For $\lambda\in\Lambda_i$ and $a\not=\pm1$ we obtain
\begin{multline*}
\widehat{\cal{L}}_z(r_{\lambda,a}\xi_{\lambda,a}+r_{\lambda,a^{-1}}\xi_{\lambda,
a^{-1}}+r_{-\lambda,a}\xi_{-\lambda,a}+r_{-\lambda,a^{-1}}\xi_{-\lambda,a^{-1}}
)\\=2\re\left(\frac{ir_{\lambda,a}\overline{r}_{\lambda,a^{-1}}}{ae^{
-2i\lambda(\eta)}-1}[\xi_{\lambda,a},\xi_{-\lambda,a^{-1}}]\right)
+2\re\left(\frac{ir_{-\lambda,a}\overline{r}_{-\lambda,a^{-1}}}{a
e^{2i\lambda(\eta)}-1}[\xi_{-\lambda,a},\xi_{\lambda,a^{-1}}]\right).
\end{multline*}
\item For $\lambda\in\Lambda_c$ and $a=1$ we obtain
\begin{multline*}
\widehat{\cal{L}}_z\bigl(r_\lambda\xi_{\lambda,1}+s_\lambda\sigma_2
(\xi_{\lambda,1})+r_{-\lambda}\xi_{-\lambda,1}+s_{-\lambda}\sigma_2
(\xi_{-\lambda,1})\bigr)\\
=2\re\left(\left(\frac{ir_\lambda\overline{s}_{-\lambda}}
{e^{-2i\lambda(\eta)}-1}-\frac{ir_{-\lambda}\overline{s}_\lambda}
{e^{2i\lambda(\eta)}-1}\right)[\xi_{\lambda,1},\xi_{-\lambda,1}]\right).
\end{multline*}
\item For $\lambda\in\Lambda_c$ and $a=-1$ we obtain
\begin{multline*}
\widehat{\cal{L}}_z\bigl(r_\lambda\xi_{\lambda,-1}+s_\lambda\sigma_2
(\xi_{\lambda,-1})+r_{-\lambda}\xi_{-\lambda,-1}+s_{-\lambda}\sigma_2
(\xi_{-\lambda,-1})\bigr)\\
=2\re\left(\left(\frac{ir_{-\lambda}\overline{s}_\lambda}
{e^{2i\lambda(\eta)}+1}-\frac{ir_\lambda\overline{s}_{-\lambda}}
{e^{-2i\lambda(\eta)}+1}\right)[\xi_{\lambda,-1},\xi_{-\lambda,-1}]\right).
\end{multline*}
\item For $\lambda\in\Lambda_c$ and $a\not=\pm1$ we obtain
\begin{multline*}
\widehat{\cal{L}}_z\bigl(r_{\lambda,a}\xi_{\lambda,a}+s_{\lambda,a}
\sigma_2(\xi_{\lambda,a})+r_{-\lambda,a^{-1}}\xi_{-\lambda,a^{-1}}+
s_{-\lambda,a^{-1}}\sigma_2(\xi_{-\lambda,a^{-1}})\bigr)\\
=2\re\left(\frac{ir_{\lambda,a}\overline{s}_{-\lambda,a^{-1}}}
{ae^{-2i\lambda(\eta)}-1}[\xi_{\lambda,a},\xi_{-\lambda,a^{-1}}]\right)
+2\re\left(\frac{ir_{-\lambda,a^{-1}}\overline{s}_{\lambda,a}}
{a^{-1}e^{2i\lambda(\eta)}-1}[\xi_{-\lambda,a^{-1}},\xi_{\lambda,a}]\right).
\end{multline*}
\end{enumerate}
\end{prop}

\begin{proof}
The proof is a straightforward application of Theorem~\ref{LeviForm:Formula}. As
illustration we will prove the first assertion. If $\lambda$ is a
real weight, we have $\sigma_2(\lambda)=\lambda$ and therefore
$\lie{u}^\mbb{C}[\lambda,a]=\lie{u}^\mbb{C}_{\lambda,a}\oplus
\lie{u}^\mbb{C}_{-\lambda,a^{-1}}$. For arbitrary numbers
$r_\lambda,r_{-\lambda}\in\mbb{C}$ we obtain
\begin{align*}
\widehat{\cal{L}}_z(r_\lambda\xi_{\lambda,a}+r_{-\lambda}
\xi_{-\lambda,a^{-1}})&=\abs{r_\lambda}^2\widehat{\cal{L}}_z
(\xi_{\lambda,a})+r_{\lambda}\overline{r}_{-\lambda}
{\cal{L}}_z(\xi_{\lambda,a},\xi_{-\lambda,a^{-1}})\\
&+\overline{r}_\lambda r_{-\lambda}
{\cal{L}}_z(\xi_{-\lambda,a^{-1}},\xi_{\lambda,a})+\abs{r_{-\lambda}}^2
\widehat{\cal{L}}_z(\xi_{-\lambda,a^{-1}})\\
&=2\re\bigl(r_\lambda\overline{r}_{-\lambda}{\cal{L}}_z(\xi_{\lambda,a},
\xi_{-\lambda,a^{-1}})\bigr)\\
&=2\re\left(\frac{ir_\lambda\overline{r}_{-\lambda}}{ae^{-2i\lambda(\eta)}-1}
[\xi_{\lambda,a},\xi_{-\lambda,a^{-1}}]\right)\\
&=2\re\left(\frac{ir_\lambda\overline{r}_{-\lambda}}{ae^{-2i\lambda(\eta)}-1}
\right)[\xi_{\lambda,a},\xi_{-\lambda,a^{-1}}],
\end{align*}
since $\sigma_2[\xi_{\lambda,a},\xi_{-\lambda,a^{-1}}]=[\xi_{\lambda,a},
\xi_{-\lambda,a^{-1}}]$ for real weights $\lambda$.
\end{proof}

\subsection{Reduction to the $(\sigma_1,\sigma_2)$--irreducible case}

In this subsection we will introduce the appropriate reduction method in order
to facilitate the determination of the Levi cone.

\begin{defn}
We say that $\lie{u}^\mbb{C}$ is $(\sigma_1,\sigma_2)$--irreducible if there is
no non-trivial ideal in $\lie{u}^\mbb{C}$ which is invariant under $\sigma_1$
and $\sigma_2$.
\end{defn}

\begin{rem}
Let $n\in A_0^c\subset U$ and $\sigma_1':=\Ad(n^{-1})\sigma_1\Ad(n)$.
Then $\sigma_1'$ is again a $\mbb{C}$--anti-linear involutive automorphism
of $\lie{u}^\mbb{C}$ commuting with $\theta$ and $\lie{u}^\mbb{C}$ is
$(\sigma_1',\sigma_2)$--irreducible if and only if it is
$(\sigma_1,\sigma_2)$--irreducible.
\end{rem}

The next lemma characterizes $(\sigma_1,\sigma_2)$--irreducibility in terms of
the set of weights $\Lambda=\Lambda(\lie{u}^\mbb{C},\lie{c}_0)$ where
$\lie{c}_0=\lie{t}_0\oplus\lie{a}_0$ is the fundamental Cartan subalgebra of
$\lie{g}_1\cap\lie{g}_2$.

\begin{lem}\label{reduction}
The Lie algebra $\lie{u}^\mbb{C}$ is $(\sigma_1,\sigma_2)$--irreducible
if and only if the root system $\Delta:=\Lambda\setminus\{0\}\subset
(i\lie{t}_0\oplus\lie{a}_0)^*$ is irreducible.
\end{lem}

\begin{proof}
Let us assume that $\lie{u}^\mbb{C}$ is $(\sigma_1,\sigma_2)$--irreducible.
If the root system $\Delta$ is not irreducible, there is a decomposition
$\Delta=\Delta_1\,\dot\cup\,\Delta_2$ into non-empty subsystems $\Delta_1$,
$\Delta_2$ such that for all $\lambda_j\in\Delta_j$ neither of
$\lambda_1\pm\lambda_2$ is a root. It follows that
\begin{equation*}
\lie{u}^\mbb{C}_j:=\lie{u}^\mbb{C}_{0,j}\oplus \bigoplus_{\lambda\in\Delta_j}
\lie{u}^\mbb{C}_\lambda,
\end{equation*}
where $\lie{u}^\mbb{C}_{0,j}:=\Span\bigl\{[\lie{u}^\mbb{C}_\lambda,
\lie{u}^\mbb{C}_{-\lambda}];\ \lambda\in\Delta_j\bigr\}$, is a non-trivial
ideal invariant under $\sigma_1$ and $\sigma_2$, which contradicts the fact
that $\lie{u}^\mbb{C}$ is $(\sigma_1,\sigma_2)$--irreducible.

In order to prove the converse, let us assume that $\lie{u}^\mbb{C}_1$ is a
non-trivial ideal in $\lie{u}^\mbb{C}$ invariant under $\sigma_1$ and
$\sigma_2$. Consequently, its orthogonal complement $\lie{u}^\mbb{C}_2$ with
respect to the Killing form $B_{\lie{u}^\mbb{C}}$ is also a non-trivial
$\sigma_1$-- and $\sigma_2$--stable ideal and $\lie{u}^\mbb{C}=\lie{u}^\mbb{C}_1
\oplus\lie{u}^\mbb{C}_2$. It is not hard to see that this decomposition induces
similar decompositions of $\lie{g}_1\cap\lie{g}_2$, $\lie{c}_0$, and hence also
of the root system $\Delta$ which contradicts the fact that $\Delta$ is
irreducible.
\end{proof}

Since the computation of the Levi form is local and the Levi form is invariant
under local biholomorphisms, it does no harm to go over to coverings. Hence, we
assume that $U^\mbb{C}$ is simply connected.

\begin{thm}
There exists an up to re-ordering unique decomposition
$\lie{u}^\mbb{C}=\lie{u}^\mbb{C}_1\oplus\dotsb\oplus\lie{u}^\mbb{C}_N$ into
$(\sigma_1,\sigma_2)$--irreducible ideals. If $U^\mbb{C}$ is simply-connected,
we have the corresponding decomposition of $U^\mbb{C}$, of the real forms $G_1$
and $G_2$, and of the orbits and their (complex) tangent spaces. This
decomposition of the complex tangent space of a closed orbit is orthogonal with
respect to its Levi form. Consequently, the Levi cone is the direct product of
the Levi cones of each factor.
\end{thm}

\begin{proof}
Since $\lie{u}^\mbb{C}$ is semi-simple, it is the direct sum of its
simple ideals, and each of these is $\theta$--invariant. Since $\sigma_1$ and
$\sigma_2$ are automorphisms of $\lie{u}^\mbb{C}$, they map simple ideals onto
simple ideals. This observation proves that $\lie{u}^\mbb{C}$ has a unique
decomposition into $(\sigma_1,\sigma_2)$--irreducible ideals. Moreover, the
simple ideals which appear in one $(\sigma_1,\sigma_2)$--irreducible ideal must
be all isomorphic.

Let $\lie{u}^\mbb{C}=\lie{u}^\mbb{C}_1\oplus\dotsb\oplus\lie{u}^\mbb{C}_N$
denote this decomposition and let $U^\mbb{C}_k$ be the subgroup of $U^\mbb{C}$
with Lie algebra $\lie{u}^\mbb{C}_k$. Since $U^\mbb{C}$ is simply-connected, we
obtain
\begin{equation*}
U^\mbb{C}\cong U^\mbb{C}_1\times\dotsb\times U^\mbb{C}_N,
\end{equation*}
and since each semi-simple normal subgroup $U^\mbb{C}_k$ is invariant under
$\sigma_1$ and $\sigma_2$, we have similar decompositions
\begin{equation*}
G_j\cong(G_j)_1\times\dotsb\times(G_j)_N
\end{equation*}
for $j=1,2$. Here, $(G_j)_k$ is the fixed point set of
$\sigma_j|_{U^\mbb{C}_k}$. It follows that the $(G_1\times G_2)$--orbits are
also direct products of their intersections with the normal subgroups
$U^\mbb{C}_k$. Since the $\lie{u}^\mbb{C}_k$ are ideals, we have
a corresponding decomposition of the set of weights into strongly orthogonal
subsystems. Finally, the computation of the Levi form in
Theorem~\ref{LeviForm:Formula} tells us that the respective parts of the complex
tangent spaces are Levi-orthogonal.
\end{proof}

\subsection{The Levi cone}

In this subsection we will determine the full Levi cone of a generic $(G_1\times
G_2)$--orbit. We assume that $U^\mbb{C}$ is simply-connected and that
$\lie{u}^\mbb{C}$ is $(\sigma_1,\sigma_2)$--irreducible.

Let $z=n\exp(i\eta)$ be a regular element contained in the standard Cartan slice
$C:=n\exp(i\lie{c})$. Since $z$ is regular, we have $e^{2i\lambda(\eta)}\not=1$
for all $(\lambda,1)\in\wt{\Lambda}$. Hence, we conclude $\lambda(\eta)\not=0$
for all $(\lambda,1)\in\wt{\Lambda}_i$. The following lemma is then a direct
consequence of Proposition~\ref{QuadrLeviForm}.

\begin{lem}\label{LeviConeGenerators}
The Levi cone ${\cal{C}}_z$ of the generic orbit $(G_1\times G_2)\cdot z$
is generated by
\begin{enumerate}[(1)]
\item $\pm[\xi_{\lambda,a},\xi_{-\lambda,a^{-1}}]$ for
$(\lambda,a)\in\wt{\Lambda}_r$,
\item $-i[\xi_{\lambda,1},\xi_{-\lambda,1}]$ for
$(\lambda,1)\in\wt{\Lambda}_i$ with $\lambda(\eta)>0$,
\item $i[\xi_{\lambda,1},\xi_{-\lambda,1}]$ for
$(\lambda,1)\in\wt{\Lambda}_i$ with $\lambda(\eta)<0$,
\item $\pm i[\xi_{\lambda,-1},\xi_{-\lambda,-1}]$ for
$(\lambda,-1)\in\wt{\Lambda}_i$,
\item $\pm\re\bigl([\xi_{\lambda,a},\xi_{-\lambda,a^{-1}}]\bigr)$ and
$\pm\im\bigl([\xi_{\lambda,a},\xi_{-\lambda,a^{-1}}]\bigr)$ for
$(\lambda,a)\in\wt{\Lambda}_i$
with $\lambda\not=\pm1$, and
\item $\pm\re\bigl([\xi_{\lambda,a},\xi_{-\lambda,a^{-1}}]\bigr)$ and
$\pm\im\bigl([\xi_{\lambda,a},\xi_{-\lambda,a^{-1}}]\bigr)$ for
$(\lambda,a)\in\wt{\Lambda}_c$.
\end{enumerate}
\end{lem}

\begin{rem}
Since we have defined the real structure on $\lie{c}^\mbb{C}$ via $\sigma_2$, we
obtain
\begin{align*}
\re\bigl([\xi_{\lambda,a},\xi_{-\lambda,a^{-1}}]\bigr)&=[\xi_{\lambda,
\lambda},\xi_{-\lambda,a^{-1}}]
+\sigma_2\bigl([\xi_{\lambda,a},\xi_{-\lambda,a^{-1}}]\bigr)\\
&=[\xi_{\lambda,a},\xi_{-\lambda,a^{-1}}]+[\xi_{\sigma_2(\lambda),
\lambda},\xi_{-\sigma_2(\lambda),\lambda^{-1}}].
\end{align*}
The imaginary part
$\im\bigl([\xi_{\lambda,a},\xi_{-\lambda,a^{-1}}]\bigr)$ can be
expressed
by an analogous formula.
\end{rem}

In order to state the main theorem we have to review some properties of real
simple Lie algebras of Hermitian type. For a more detailed exposition of these
topics we refer the reader to~\cite{HiNe} and~\cite{Ne2}.

Recall that a simple real Lie algebra $\lie{g}=\lie{k}\oplus\lie{p}$ is said to
be of Hermitian type if the center of $\lie{k}$ is non-trivial. This condition
implies that a maximal torus $\lie{t}\subset\lie{k}$ is a Cartan subalgebra of
$\lie{g}$. Then every root $\alpha$ in
$\Delta=\Delta(\lie{g}^\mbb{C},\lie{t})$ is imaginary, and either
$\lie{g}^\mbb{C}_\alpha\subset\lie{k}^\mbb{C}$ or $\lie{g}^\mbb{C}_\alpha
\subset\lie{p}^\mbb{C}$ holds. In the first case we call $\alpha$ a compact
root, while in the second case $\alpha$ is said to be non-compact. We write
$\Delta_\lie{k}$ and $\Delta_\lie{p}$ for the sets of compact and non-compact
roots, respectively. Since $\lie{g}$ is Hermitian, the root system $\Delta$
possesses a good ordering, i.\,e.\ there is a choice of the set $\Delta^+$ of
positive roots such that each positive non-compact root is larger than every
compact root. This is equivalent to the fact that the set $\Delta_\lie{p}^+$ is
invariant under the Weyl group $W(\Delta_\lie{k})$. Therefore there are two
natural $W(\Delta_\lie{k})$--invariant cones $C_{\min}\subset C_{\max}$, where
$C_{\min}$ is the closed convex cone generated by
\begin{equation*}
\bigl\{-i\bigl[\xi_\alpha,\sigma(\xi_\alpha)\bigr];\
\xi_\alpha\in\lie{g}^\mbb{C} _\alpha, \alpha\in\Delta_\lie{p}^+\bigr\}
\subset\lie{t}
\end{equation*}
and
\begin{equation*}
C_{\max}:=\bigl\{\eta\in\lie{t};\ i\alpha(\eta)\geq0\text{ for all
$\alpha\in\Delta_\lie{p}^+$}\bigr\}.
\end{equation*}
Let $C_{\max}^0$ be the interior of $C_{\max}$. Then the open subset
$G\exp(iC_{\max}^0)G\subset U^\mbb{C}$ is closed under multiplication
and hence a semi-group, called the open complex Ol'shanski{\u\i} semi-group.

\begin{thm}\label{LeviCone}
Let $\lie{u}^\mbb{C}$ be $(\sigma_1,\sigma_2)$--irreducible and let $(G_1\times
G_2)\cdot z$ be a generic orbit where $z=n\exp(i\eta)$ lies in the standard
Cartan slice $C:=n\exp(i\lie{c})$.
\begin{enumerate}[(1)]
\item If the standard Cartan subset $\lie{c}$ is non-compact, then
${\cal{C}}_z=\lie{c}$ holds.
\item If $\lie{c}$ is compact and if $a\not=1$ for some
$(\lambda,a)\in\wt{\Lambda}$, then we have ${\cal{C}}_z=\lie{c}$.
\item If $\lie{c}$ is compact and if $a=1$ for all generalized weights,
then $\sigma_1=\sigma_2$ holds and there are the following cases.
\begin{enumerate}[(i)]
\item If $\lie{g}_1=\lie{g}_2=:\lie{g}$ is of Hermitian type and if $\eta$ lies
in $C_{\max}$, then the Levi cone ${\cal{C}}_z$ is isomorphic to the dual of the
positive Weyl chamber defined by $\Lambda^+$. In particular, the Levi cone is
pointed.
\item If $\lie{g}$ is of Hermitian type and $\eta\notin C_{\max}$, then
${\cal{C}}_z=\lie{c}$.
\item If $\lie{g}$ is not of Hermitian type, then ${\cal{C}}_z=\lie{c}$.
\end{enumerate}
\end{enumerate}
\end{thm}

\begin{rem}
The reader will note that the statement of Theorem~\ref{LeviCone} differs also
for the case
$\sigma_1=\sigma_2$ from the corresponding Theorem~5.3 in~\cite{FeGe}.
Indeed, as L.~Geatti has kindly pointed out, the formulation of the third part
of Theorem~5.3
in~\cite{FeGe} is not correct. The correct statement in Theorem~\ref{LeviCone}
and its proof in
the case $\sigma_1=\sigma_2$ are due to an unpublished erratum written by
L.~Geatti.
\end{rem}

It will turn out to be convenient to express the generators of the Levi cone in
terms of the coroots $\eta_\lambda\in i\lie{t}\oplus\lie{a}$. Therefore we will
identify $\lie{c}=\lie{t}\oplus\lie{a}$ with $i\lie{t}\oplus\lie{a}$
via the map $(\eta_1,\eta_2)\mapsto(i\eta_1,\eta_2)$. By abuse of notation, we
denote the image of the Levi cone under this map again by ${\cal{C}}_z\subset
i\lie{t}\oplus\lie{a}$.
According to Lemma~\ref{WeightVectorProperties} we have
\begin{equation*}
[\xi_{\lambda,a},\xi_{-\lambda,a^{-1}}]=
B_{\lie{u}^\mbb{C}}(\xi_{\lambda,a},\xi_{-\lambda,a^{-1}}
)\eta_\lambda\in\mbb{C}\eta_\lambda.
\end{equation*}
Hence, we can normalize the $\xi_{\lambda,a}$ such that
$[\xi_{\lambda,a},\xi_{-\lambda,a^{-1}}]=\eta_\lambda$
holds for all $\lambda\in\Lambda^+\setminus\Lambda_i$ and
$[\xi_{\lambda,a},\xi_{-\lambda,a^{-1}}]=\pm\eta_\lambda$
holds for $\lambda\in\Lambda_i^+$ depending on the sign of
$B_{\lie{u}^\mbb{C}}(\xi_{\lambda,a},\xi_{-\lambda,a^{-1}})$.

\begin{rem}
In the case where $\lie{g}_1=\lie{g}_2=:\lie{g}$ and $\lie{t}$ is a compact
Cartan subalgebra of $\lie{g}$, we obtain after the above normalization
\begin{equation*}
[\xi_\alpha,\xi_{-\alpha}]=
\begin{cases}
\eta_\alpha & \text{ for $\alpha\in\Delta_\lie{p}^+$}\\
-\eta_\alpha & \text{ for $\alpha\in\Delta_\lie{k}^+$}
\end{cases},
\end{equation*}
since the real part of $B_{\lie{u}^\mbb{C}}$ is positive definite on $\lie{p}$
and negative definite on $\lie{k}$.
\end{rem}

\begin{proof}[Proof of Theorem~\ref{LeviCone}]

{\bf(1)} Let $\lie{c}$ be non-compact. Since $\Lambda\setminus\{0\}$ satisfies
the axioms for an abstract root system, we may choose a set
$\Pi\subset\Lambda^+$ of simple weights. By
Lemma~\ref{LeviConeGenerators} we know that $\pm\eta_\lambda$ lies in
${\cal{C}}_z$ for $\lambda\in\Lambda\setminus\Lambda_i$, and we have to show
that $\pm\eta_\lambda\in{\cal{C}}_z$ holds for all $\lambda\in\Lambda$. It is
enough to prove this fact for all $\eta_\lambda$ with
$\lambda\in\Pi_i:=\Pi\cap\Lambda_i$.

If $\lambda,\mu\in\Pi_i$ with $\lambda+\mu\in\Lambda$ are given, then
$\lambda+\mu\in\Lambda_i^+$ holds. Since $\lie{c}$ is non-compact, this
observation implies that there exists an element $\mu\in\Pi\setminus\Pi_i$.
Let $\lambda\in\Pi_i$ be arbitrary (if $\Pi_i=\emptyset$, the proof is
finished). Since $\Lambda\setminus\{0\}$ is irreducible by
Lemma~\ref{reduction}, its Dynkin diagram is connected and hence we find a
sequence $\lambda=\lambda_1,\dotsc,\lambda_N=\mu$ of simple roots which are
adjacent in the Dynkin diagram. Consequently, we obtain
$\lambda_j+\dotsb+\lambda_N\in\Lambda\setminus\Lambda_i$ for all $0\leq j\leq
N-1$. This implies
\begin{equation*}
\pm\eta_{\lambda_j+\dotsb+\lambda_N}=\pm(\eta_{\lambda_j+\dotsb+\lambda_{N-1}}
+\eta_{\lambda_N})\in{\cal{C}}_z
\end{equation*}
for all $0\leq j\leq N-1$. Since $\pm\eta_{\lambda_N}$ lies in ${\cal{C}}_z$, we
conclude $\pm\eta_{\lambda_j+\dotsb+\lambda_{N-1}}\in{\cal{C}}_z$ for all $j$.
Iterating this argument we finally arrive at $\pm\eta_{\lambda}\in{\cal{C}}_z$
which was to be shown.

{\bf(2)} Let us assume that $\lie{c}$ is compact and that there exists
$(\lambda,a)\in\wt{\Lambda}$ with $a\not=1$. In this case we
have $\Lambda=\Lambda_i$ and $\pm\eta_\lambda\in{\cal{C}}_z$ for all $\lambda$
such that there exists $a\not=1$ with $(\lambda,a)\in\wt{\Lambda}$.
If there are two weights $\lambda_1,\lambda_2\in\Lambda$ such that
$(\lambda_j,a)\in\wt{\Lambda}$ implies $a=1$ for $j=1,2$ and
such that $\lambda_1+\lambda_2$ is again a weight, then we conclude from
Proposition~\ref{WeightSpaceDecomposition} that
\begin{equation*}
(\lambda_1+\lambda_2,a)\in\wt{\Lambda}\Longrightarrow a=1
\end{equation*}
holds. Consequently, each set $\Pi\subset\Lambda^+$ of simple roots must contain
a root $\mu$ with $\pm\eta_\mu\in{\cal{C}}_z$. Now the claim follows from the
same argument as above.

{\bf(3)} Let $\lie{c}=\lie{t}$ be a compact Cartan subalgebra of
$\lie{g}_2\cap\Ad(n^{-1})\lie{g}_1$ such that $a=1$ holds for each extended
weight $(\lambda,a)\in\wt{\Lambda}$. It is enough to prove that this assumption
implies $\lie{g}_1=\lie{g}_2$ since then the claim follows from~\cite{FeGe} and
Geatti's erratum.

The proof of $\lie{g}_1=\lie{g}_2$ relies on the comparison of the weight space
decompositions
\begin{equation*}
\lie{u}^\mbb{C}=\lie{t}^\mbb{C}\oplus\bigoplus_{\lambda\in\Lambda\setminus
\{0\}}\lie{u}^\mbb{C}_{\lambda,1}\quad\text{and}\quad
\bigl(\lie{g}_2\cap\Ad(n^{-1})\lie{g}_1)^\mbb{C}=\lie{t}^\mbb{C}\oplus
\bigoplus_{\lambda\in\Lambda'}\lie{u}^\mbb{C}_{\lambda,1},
\end{equation*}
where $\Lambda'$ denotes the set of non-zero weights $\lambda$ for which
$\lie{u}^\mbb{C}_{\lambda,1}$ is contained in $\bigl(\lie{g}_2\cap\Ad(n^{-1})
\lie{g}_1\bigr)^\mbb{C}$. Note that this is well-defined since
$\dim\lie{u}^\mbb{C}_{\lambda,1}=1$ by
Proposition~\ref{WeightSpaceDecomposition}. Since the weight space
decomposition is in both cases defined with respect to $\lie{t}$, a basis of
$\Lambda'$ has to be a basis of $\Lambda\setminus\{0\}$, too. Since the root
system $\Lambda\setminus\{0\}$ is reduced by
Proposition~\ref{WeightSpaceDecomposition}, we conclude
\begin{equation*}
\lie{u}^\mbb{C}=\bigl(\lie{g}_2\cap\Ad(n^{-1})\lie{g}_1\bigr)^\mbb{C},
\end{equation*}
and hence, that $\lie{g}_2\cap\Ad(n^{-1})\lie{g}_1$ is a real form of
$\lie{u}^\mbb{C}$. For dimensional reasons this implies
\begin{equation*}
\lie{g}_2=\lie{g}_2\cap\Ad(n^{-1})\lie{g}_1=\Ad(n^{-1})\lie{g}_1,
\end{equation*}
i.\,e.\ $\sigma_2=\Ad(n^{-1})\sigma_1\Ad(n)$.

By the definition of a standard Cartan subset, the fundamental Cartan
subalgebra $\lie{c}_0\subset\lie{g}_1\cap\lie{g}_2$ has the same dimension as
$\lie{t}$. Therefore, we obtain
\begin{equation*}
\rk(\lie{g}_1\cap\lie{g}_2)=\dim\lie{c}_0=\dim\lie{t}=
\rk\bigl(\lie{g}_2\cap\Ad(n^{-1})\lie{g}_1\bigr)^\mbb{C}=\rk\lie{u}^\mbb{C},
\end{equation*}
which implies in the same way as above that $\lie{g}_1\cap\lie{g}_2$ is a real
form of $\lie{u}^\mbb{C}$ and hence that $\lie{g}_1=\lie{g}_2$ holds.
\end{proof}

\section{Applications}

\subsection{The criterion of Fels and Geatti}

We restate Corollary~5.6 from~\cite{FeGe} whose proof relies
on Theorem~\ref{BoggessPolking}.

\begin{thm}[Fels, Geatti]\label{Criterion}
Let $Z$ be a complex manifold on which the Lie group $G$ acts by holomorphic
transformations. Let the orbit $M_z=G\cdot z$ be a generic CR submanifold such
that $\cal{C}_z=T_zM_z/H_zM_z$ holds. Then there exists no $G$--invariant Stein
domain in $Z$ which contains $M_z$ in its boundary. Furthermore, there is no
non-constant $G$--invariant plurisubharmonic function which is defined in a
neighborhood of $M_z$.
\end{thm}

Theorem~\ref{Criterion} gives a necessary condition for an invariant domain with
a generic orbit in its boundary to be Stein. In our situation we obtain the
following result.

\begin{thm}
Let $C=n\exp(i\lie{c})$ be a standard Cartan subset and let $\Omega$ be a
connected component of the open set $G_1(C\cap U^\mbb{C}_{sr})G_2$. Then
$\Omega$ does not contain any proper $(G_1\times G_2)$--invariant Stein
subdomain unless $\lie{c}$ is compact and $\tau_n$ has $a=1$ as only
eigenvalue. In this case $G_1=G_2$ must be of Hermitian type and $\Omega$ is a
translate of the open Ol'shanski{\u\i} semi-group in $U^\mbb{C}$.
\end{thm}

Consequently, we see that in the case $G_1\not=G_2$ there are no invariant
Stein subdomains in $U^\mbb{C}$ in whose boundary a generic orbit lies. The
reader should note that there are only finitely many $(G_1\times
G_2)$--invariant domains whose boundaries consist entirely of non-generic
orbits.

\subsection{$q$--pseudo-convex functions and $q$--completeness}

In this subsection we review quickly the notions of $q$--pseudo-convex
functions and $q$--complete complex manifolds. Let $\Omega$ be a domain in a
complex manifold $Z$. We call a smooth function on $\Omega$ strictly
$q$--pseudo-convex if its Levi form has at least $n-q$ positive eigenvalues,
$n:=\dim_\mbb{C}\Omega$, at each point of $\Omega$. Hence, a strictly
$0$--pseudo-convex function is the same as a strictly plurisubharmonic function.
If $\Omega$ admits a strictly $q$--pseudo-convex exhaustion function, we say
that $\Omega$ is $q$--complete. The solution of the Levi problem implies that a
domain is Stein if and only if it is $0$--complete. For more properties of
$q$--complete complex spaces we refer the reader to~\cite{AnGr}.

\begin{rem}
A standard argument of complex analysis (compare Corollary~XIII.5.4
in~\cite{Ne2} for the case $q=0$ and~\cite{Dem} for the generalization to $q>0$)
shows that a domain $\Omega$ in a Stein manifold $Z$ is $q$--complete if and
only if there exists a strictly $q$--pseudo-convex function $\Phi$ on $\Omega$
with the property $\Phi(z_n)\to\infty$ whenever $z_n\to z\in\partial\Omega$.
\end{rem}

A domain $\Omega\subset Z$ with smooth boundary is called Levi--$q$--convex, if
the boundary $\partial\Omega$ can locally be defined by a function whose Levi
form has at most $q$ negative eigenvalues when restricted to the complex tangent
space at any point of $\partial\Omega$.

By a theorem of Docquier and Grauert (\cite{DocGr}) a domain $\Omega$ with
smooth boundary in a Stein manifold $Z$ is Stein if and only if it is
Levi--$0$--convex. In~\cite{EVS} this result is generalized to arbitrary $q$.

\begin{thm}[Oka, Docquier-Grauert, Eastwood-Suria]\label{Levi-q-convex}
Let $Z$ be a Stein manifold and let $\Omega\subset Z$ be a domain with smooth
boundary. Then $\Omega$ is strictly $q$--complete if and only if $\Omega$ is
Levi--$q$--convex.
\end{thm}

\subsection{The rank one case}

If the closed orbit $M_z=(G_1\times G_2)\cdot z$ is a hypersurface in
$U^\mbb{C}$, its intrinsic Levi form coincides with the classical Levi form of
that hypersurface, and hence the signature of $\widehat{\mathcal{L}}_z$ is
defined. According to Theorem~\ref{Levi-q-convex} this signature encodes
information about complex-analytic properties of the domains bounded by $M_z$.
In this subsection we will use Matsuki's classification of pairs of
involutive automorphisms of simply-connected compact Lie groups in order to
classify the triples $(U^\mbb{C},G_1,G_2)$ where $U^\mbb{C}$ is
simply-connected and the generic $(G_1\times G_2)$--orbit is a hypersurface.
Moreover, we will determine the signature of the Levi form of each generic
hypersurface orbit.

In~\cite{Ma2} pairs of involutive automorphisms of simply-connected semi-simple
compact Lie groups are classified under the following notion of equivalence.

\begin{defn}
Let $U$ be a simply-connected semi-simple compact Lie group. Two pairs of
involutive automorphisms $(\sigma_1,\sigma_2)$ and $(\sigma_1',\sigma_2')$ are
called equivalent if there exist an automorphism $\varphi\in\Aut(U)$ and an
element $u\in U$ such that
\begin{equation*}
\sigma_1'=\varphi\sigma_1\varphi^{-1}\quad\text{and}\quad
\sigma_2'=\Int(u)\varphi\sigma_2\varphi^{-1}\Int(u)^{-1}
\end{equation*}
hold.
\end{defn}

Since in our case the involutive automorphisms $\sigma_1,\sigma_2\colon
U^\mbb{C}\to U^\mbb{C}$ commute with $\theta$ and are anti-holomorphic, they are
completely determined by their restrictions to $U$. Therefore, we may apply the
classification result from~\cite{Ma2}.

\begin{thm}[Matsuki]\label{Classification}
Let $U^\mbb{C}$ be simply-connected. If the generic $(G_1\times G_2)$--orbit is
a hypersurface in $U^\mbb{C}$, then $U^\mbb{C}$ is of the form
\begin{equation*}
U^\mbb{C}=\underbrace{S\times\dotsb\times S}_{\text{$k$ times}},
\end{equation*}
where $S$ is a $\theta$--stable normal subgroup of $U^\mbb{C}$ either
isomorphic to ${\rm{SL}}(2,\mbb{C})$ or ${\rm{SL}}(3,\mbb{C})$. Let $\sigma$
and $\tau$ be anti-holomorphic involutive automorphisms of $S$ commuting with
$\theta|_S$. If $k$ is odd,
then we have
\begin{align*}
\sigma_1(g_1,\dotsc,g_k)&=\bigl(\sigma(g_1),\theta(g_3),\theta(g_2),
\dotsc,\theta(g_k),\theta(g_{k-1})\bigr)\\
\sigma_2(g_1,\dotsc,g_k)&=\bigl(\theta(g_2),\theta(g_1),\dotsc,
\theta(g_{k-1}),\theta(g_{k-2}),\tau(g_k)\bigr),
\end{align*}
and if $k$ is even, then
\begin{align*}
\sigma_1(g_1,\dotsc,g_k)&=\bigl(\sigma(g_1),\theta(g_3),\theta(g_2),
\dotsc,\theta(g_{k-1}),\theta(g_{k-2}),\tau(g_k)\bigr)\\
\sigma_2(g_1,\dotsc,g_k)&=\bigl(\theta(g_2),\theta(g_1),\dotsc,
\theta(g_k),\theta(g_{k-1})\bigr)
\end{align*}
holds. If $S={\rm{SL}}(2,\mbb{C})$, then the pair $(\sigma,\tau)$ is equivalent
to one of $\bigl\{(\sigma_{1,1},\sigma_{1,1}),(\sigma_{1,1},\theta),
(\theta,\theta)\bigr\}$, where $\sigma_{1,1}$ is the involution defining the
non-compact real form ${\rm{SU}}(1,1)$ of ${\rm{SL}}(2,\mbb{C})$. If
$S={\rm{SL}}(3,\mbb{C})$, then the only possibility for $(\sigma,\tau)$ up to
equivalence is $\bigl(\sigma(g),\tau(g)\bigr)=\bigl(\overline{g},I_{2,1}
\theta(g)I_{2,1}\bigr)$ with $I_{2,1}:=\left(\begin{smallmatrix}1&0&0\\0&1&0\\
0&0&-1\end{smallmatrix}\right)$.
\end{thm}

\begin{proof}
Since the semi-simple complex Lie group $U^\mbb{C}$ is assumed to be
simply-connected, we can identify the automorphism group $\Aut(U^\mbb{C})$ with
$\Aut(\lie{u}^\mbb{C})$. By Proposition~2.2 in~\cite{Ma2} there exists a
$\theta$--invariant decomposition
\begin{equation*}
\lie{u}^\mbb{C}=\lie{u}^\mbb{C}_1\oplus\dotsb\oplus\lie{u}^\mbb{C}_N
\end{equation*}
into $\sigma_1$-- and $\sigma_2$--invariant semi-simple ideals
$\lie{u}^\mbb{C}_j$. Moreover, each $\lie{u}^\mbb{C}_j$ is of the form
\begin{equation*}
\lie{u}^\mbb{C}_j=\underbrace{\lie{s}_j\oplus\dotsb\oplus\lie{s}_j}_{\text{$k_j$
times}},
\end{equation*}
where $\lie{s}_j$ is a $\theta$--stable simple ideal in $\lie{u}^\mbb{C}_j$,
such that the restriction of the pair $(\sigma_1,\sigma_2)$ (or
$(\sigma_2,\sigma_1)$) to $\lie{u}^\mbb{C}_j$ is equivalent to one of the
following three types:
\begin{enumerate}[(1)]
\item The number $k_j$ is even and
\begin{align*}
\sigma_1(\xi_1,\dotsc,\xi_k)&=\bigl(\varphi(\xi_k),\theta(\xi_3),\theta(\xi_2),
\dotsc,\theta(\xi_{k-1}),\theta(\xi_{k-2}),\varphi^{-1}(\xi_1)\bigr)\\
\sigma_2(\xi_1,\dotsc,\xi_k)&=\bigl(\theta(\xi_2),\theta(\xi_1),\dotsc,
\theta(\xi_k),\theta(\xi_{k-1})\bigr),
\end{align*}
for some $\mbb{C}$--anti-linear automorphism $\varphi$ of $\lie{s}_j$ commuting
with $\theta|_{\lie{s}_j}$.
\item The number $k$ is even and
\begin{align*}
\sigma_1(\xi_1,\dotsc,\xi_k)&=\bigl(\sigma(\xi_1),\theta(\xi_3),\theta(\xi_2),
\dotsc,\theta(\xi_{k-1}),\theta(\xi_{k-2}),\tau(\xi_k)\bigr)\\
\sigma_2(\xi_1,\dotsc,\xi_k)&=\bigl(\theta(\xi_2),\theta(\xi_1),\dotsc,
\theta(\xi_k),\theta(\xi_{k-1})\bigr),
\end{align*}
where $\sigma$ and $\tau$ are $\mbb{C}$--anti-linear involutive automorphisms of
$\lie{s}_j$ commuting with $\theta|_{\lie{s}_j}$.
\item The number $k$ is odd and
\begin{align*}
\sigma_1(\xi_1,\dotsc,\xi_k)&=\bigl(\sigma(\xi_1),\theta(\xi_3),\theta(\xi_2),
\dotsc,\theta(\xi_k),\theta(\xi_{k-1})\bigr)\\
\sigma_2(\xi_1,\dotsc,\xi_k)&=\bigl(\theta(\xi_2),\theta(\xi_1),\dotsc,
\theta(\xi_{k-1}),\theta(\xi_{k-2}),\tau(\xi_k)\bigr),
\end{align*}
where $\sigma$ and $\tau$ are $\mbb{C}$--anti-linear involutive automorphisms of
$\lie{s}_j$ commuting with $\theta|_{\lie{s}_j}$.
\end{enumerate}
The condition that the generic $(G_1\times G_2)$--orbit is a hypersurface is
equivalent to $\rk(\lie{g}_1\cap\lie{g}_2)=1$. In particular, this condition
implies that $\lie{u}^\mbb{C}$ is $(\sigma_1,\sigma_2)$--irreducible.

If $\lie{u}^\mbb{C}$ is of the first type, one checks directly that
\begin{equation*}
\lie{g}_1\cap\lie{g}_2\cong(\lie{s}_j^{\theta\varphi})^\mbb{R}
\end{equation*}
holds. Consequently $\rk(\lie{g}_1\cap\lie{g}_2)$ is even and in particular
larger than $1$. This excludes the first type.

Let $\lie{u}^\mbb{C}$ be of the second or the third type. Again it is not hard
to see that
\begin{equation*}
\lie{g}_1\cap\lie{g}_2\cong\lie{s}_j^\sigma\cap\lie{s}_j^\tau
\end{equation*}
holds. It follows that the simple complex Lie algebra $\lie{s}_j$ contains the
complex subalgebra $(\lie{g}_1\cap\lie{g}_2)^\mbb{C}$ which is given as the set
of fixed points of the $\mbb{C}$--linear semi-simple automorphism $\sigma\tau$.
Using the classification of semi-simple automorphisms of simple complex Lie
algebras we see that the only possibilities for $\lie{s}$ are $\lie{sl}(2,
\mbb{C})$ and $\lie{sl}(3,\mbb{C})$. Then the claim follows from
Proposition~2.1 in~\cite{Ma2} where the pairs of involutions on the
classical Lie algebras are classified up to equivalence.
\end{proof}

\begin{rem}
Let $U^\mbb{C}=S\times\dotsb\times S$ ($k$ times) with $S={\rm{SL}}(2,\mbb{C})$,
and let us consider the involutions $\sigma_1$ and $\sigma_2$ on $U^\mbb{C}$
corresponding to $(\sigma_{1,1},\theta)$ in the way described in
Theorem~\ref{Classification}. In this case we see that
$\lie{g}_1\cap\lie{g}_2=:\lie{t}_0\cong\lie{so}(2,\mbb{R})$ is one-dimensional
and compact. Hence, the fundamental Cartan subset $C_0=\exp(i\lie{t}_0)$ is an
exact slice for the $(G_1\times G_2)$--action on $U^\mbb{C}$, i.\,e.\ every
$(G_1\times G_2)$--orbit intersects $C_0$ in exactly one point. In particular,
we conclude that each element $z\in U^\mbb{C}$ is strongly regular and that
$G_1\times G_2$ acts properly on $U^\mbb{C}$.
\end{rem}

In the following let us consider a point $z\in U^\mbb{C}$ such that the orbit
$M_z=(G_1\times G_2)\cdot z$ is a closed hypersurface in $U^\mbb{C}$. Without
loss of generality we take $z$ to be of the form $z=\exp(i\eta)\in C$ for some
standard Cartan subset $C=\exp(i\lie{c})$. Because of $\rk(\lie{g}_1\cap
\lie{g}_2)=1$ the Cartan subalgebra $\lie{c}\subset\lie{g}_1\cap\lie{g}_2$ is
one-dimensional and hence either $\lie{c}=\lie{t}$ is a maximal torus in
$\lie{k}_1\cap\lie{k}_2$ or $\lie{c}=\lie{a}$ is a maximal Abelian subspace of
$\lie{p}_1\cap\lie{p}_2$. Let
\begin{equation*}
\lie{u}^\mbb{C}=\bigoplus_{(\lambda,a)\in\wt{\Lambda}}
\lie{u}^\mbb{C}_{\lambda,a}
\end{equation*}
be the extended weight space decomposition of $\lie{u}^\mbb{C}$ with respect to
$\lie{c}$.

Let us first assume that $\lie{c}=\lie{a}$ is non-compact. In this case every
weight is real and we conclude from Theorem~\ref{LeviForm:Formula} that the only
non-zero contributions to the Levi form of $M_z$ stem from the terms
\begin{equation*}
\mathcal{L}_z(\xi_{\lambda,a},\xi_{-\lambda,a^{-1}})=
\frac{i}{ae^{-2i\lambda(\eta)}-1}[\xi_{\lambda,a},\xi_{-\lambda,a^{-1}}],
\end{equation*}
where $\lambda\not=0$ and $\xi_{\lambda,a}$ is a non-zero element in
$\lie{u}^\mbb{C}_{\lambda,a}$ with
$\sigma_2(\xi_{\lambda,a})=\xi_{\lambda,a}$. Consequently, the restriction of
the Levi form to
\begin{equation*}
\lie{u}^\mbb{C}[\lambda,a]=
\lie{u}^\mbb{C}_{\lambda,a}\oplus\lie{u}^\mbb{C}_{-\lambda,a^{-1}},\quad
\lambda\in\Lambda^+,
\end{equation*}
has with respect to the bases $(\xi_{\lambda,a},\xi_{-\lambda,a^{-1}})$ of
$\lie{u}^\mbb{C}[\lambda,a]$ and
$[\xi_{\lambda,a},\xi_{-\lambda,a^{-1}}]$ of $\lie{a}$ the matrix
\begin{equation*}
\begin{pmatrix}
0&\frac{i}{ae^{-2i\lambda(\eta)}-1}\\
-\frac{i}{a^{-1}e^{2i\lambda(\eta)}-1}&0
\end{pmatrix},
\end{equation*}
which has the eigenvalues $\pm\frac{1}{\abs{ae^{-2i\lambda(\eta)}-1}}$.
Hence, we obtain a pair of one positive and one negative eigenvalue of the Levi
form on $\lie{u}^\mbb{C}[\lambda,a]$ for each $\lambda\in\Lambda^+$.

If $\lie{c}=\lie{t}$ is compact, each weight is imaginary and we have to
handle the cases $a=\pm1$ and $a\not=\pm1$ separately. If $a=\pm1$, then
$a=a^{-1}$ and consequently
\begin{equation*}
\lie{u}^\mbb{C}[\lambda,a]=\lie{u}^\mbb{C}_{\lambda,a}\oplus
\lie{u}^\mbb{C}_{-\lambda,a},\quad\lambda\in\Lambda^+,
\end{equation*}
holds. As basis of $\lie{u}^\mbb{C}[\lambda,a]$ we choose
$(\xi_{\lambda,a},\xi_{-\lambda,a})$ with
$\sigma_2(\xi_{\lambda,a})=\xi_{-\lambda,a}$, and as basis of $\lie{t}$ we take
$i[\xi_{\lambda,a},\xi_{-\lambda,a}]$. Then the Levi form has with respect to
these bases the matrix
\begin{equation*}
\begin{pmatrix}
\frac{1}{ae^{-2i\lambda(\eta)}-1}&0\\
0&-\frac{1}{ae^{2i\lambda(\eta)}-1}
\end{pmatrix}.
\end{equation*}
If $a=1$, then both eigenvalues have the same sign, and if $a=-1$, then the
eigenvalues have different sign.

For $a\not=\pm1$ we have
\begin{equation*}
\lie{u}^\mbb{C}[\lambda,a]=\lie{u}^\mbb{C}_{\lambda,a}\oplus
\lie{u}^\mbb{C}_{-\lambda,a^{-1}}\oplus\lie{u}^\mbb{C}_{\lambda,a^{-1}}\oplus
\lie{u}^\mbb{C}_{-\lambda,a}
\end{equation*}
and take $(\xi_{\lambda,a},\xi_{-\lambda,a^{-1}},\xi_{\lambda,a^{-1}},
\xi_{-\lambda,a})$ as a basis of $\lie{u}^\mbb{C}[\lambda,a]$. Under the
assumption $B_{\lie{u}^\mbb{C}}(\xi_{\lambda,a},\xi_{-\lambda,a^{-1}})=
B_{\lie{u}^\mbb{C}}(\xi_{\lambda,a^{-1}},\xi_{-\lambda,a})$ we obtain
$i[\xi_{\lambda,a},\xi_{-\lambda,a^{-1}}]=i[\xi_{\lambda,a^{-1}},
\xi_{-\lambda,a}]$ which we take as a basis of $\lie{t}$. With respect to these
bases the restriction of the Levi form has the matrix
\begin{equation*}
\begin{pmatrix}
0&\frac{1}{ae^{-2i\lambda(\eta)}-1}&0&0\\
\frac{1}{a^{-1}e^{-2i\lambda(\eta)}-1}&0&0&0\\
0&0&0&\frac{1}{a^{-1}e^{-2i\lambda(\eta)}-1}\\
0&0&\frac{1}{ae^{-2i\lambda(\eta)}-1}&0
\end{pmatrix},
\end{equation*}
whose eigenvalues are given by $\pm\frac{1}{\abs{ae^{-2i\lambda(\eta)}-1}}$
and $\pm\frac{1}{\abs{a^{-1}e^{-2i\lambda(\eta)}-1}}$.

We summarize these results in the following

\begin{thm}\label{Signature}
If $\lie{c}=\lie{a}$ is non-compact, then each generic orbit $M_z$ with
$z\in\exp(i\lie{a})$ is Levi--$q$--convex with
\begin{equation*}
q=\#\Lambda^+.
\end{equation*}
If $\lie{c}=\lie{t}$ is compact, let us choose an ordering on the set of weights
such that $\lambda(\eta)<0$ for all $\lambda\in\Lambda^+$ and $z=\exp(i\eta)$.
Then each generic orbit $M_z$ with $z\in\exp(i\lie{t})$ is Levi--$q$--convex
with
\begin{equation*}
q=\#\bigl\{(\lambda,-1)\in\wt{\Lambda};\ \lambda\in\Lambda^+\bigr\}+\#\bigl\{
(\lambda,a)\in\wt{\Lambda};\ a\not=\pm1\bigr\}
\end{equation*}
Moreover, this numbers for $q$ are sharp, i.\,e.\ $M_z$ is not
Levi--$q'$--convex for any $q'<q$.
\end{thm}

\begin{proof}
The only claim which is left to show is the fact that the multiplicity of the
eigenvalue $0$ of $\mathcal{L}_z$ is given by $\rk(\lie{u}^\mbb{C})-1$.

We conclude from Theorem~\ref{LeviForm:Formula} that the nullspace of the Levi
form of a generic orbit coincides with 
\begin{equation*}
\bigoplus_{(0,a)\in\wt{\Lambda}:a\not=1}\lie{u}^\mbb{C}_{0,a}
\end{equation*}
According to Lemma~5.1 in~\cite{Ma2} the subalgebra $\lie{u}^\mbb{C}_0$ is a
Cartan subalgebra of $\lie{u}^\mbb{C}$. Since we assume that the generic $
(G_1\times G_2)$--orbit is a hypersurface, we conclude $\dim_\mbb{C}
\lie{u}^\mbb{C}_{0,1}=1$ which finishes the proof.
\end{proof}

Theorem~\ref{Levi-q-convex} and Theorem~\ref{Signature} yield the following
result.

\begin{thm}
Let $M_z$ be a closed hypersurface orbit where $z\in C=\exp(i\lie{c})$ and let
$\Omega$ be an invariant domain with $\partial\Omega=M_z$.
\begin{enumerate}[(1)]
\item If $\lie{c}=\lie{a}$ is non-compact, then $\Omega$ or $U^\mbb{C}\setminus
\overline{\Omega}$ is strictly $q$--complete with
\begin{equation*}
q=\#\Lambda^+,
\end{equation*}
and this $q$ is optimal.
\item If $\lie{c}=\lie{t}$ is compact, then $\Omega$ or $U^\mbb{C}\setminus
\overline{\Omega}$ is strictly $q$--complete with
\begin{equation*}
q=\#\bigl\{(\lambda,-1)\in\wt{\Lambda};\ \lambda\in\Lambda^+\bigr\}+\#\bigl\{
(\lambda,a)\in\wt{\Lambda};\ a\not=\pm1\bigr\}
\end{equation*}
and this $q$ is optimal.
\end{enumerate}
\end{thm}

\subsection{The Levi form of invariant functions and $q$--complete domains}

Let us assume in this subsection that the intersection $\lie{g}_1\cap\lie{g}_2$
contains a compact Cartan subalgebra $\lie{t}$. By Matsuki's result a generic
orbit $M_z=(G_1\times G_2)\cdot z$ intersects the corresponding standard Cartan
subset $C=\exp(i\lie{t})$ in an orbit of the group $W_{K_1\times K_2}(C)$.

\begin{rem}
It can be shown that there exists a group isomorphism from
\begin{equation*}
W:=W_{K_1\cap K_2}(\lie{t}):={\mathcal{N}}_{K_1\cap K_2}(\lie{t})/
\mathcal{Z}_{K_1\cap K_2}(\lie{t})
\end{equation*}
onto $W_{K_1\times K_2}(C)$ such that the diffeomorphism $\lie{t}\to C$,
$\eta\mapsto\exp(i\eta)$, intertwines the $W$--action on $\lie{t}$ with the
$W_{K_1\times K_2}(C)$--action on $C$ modulo this isomorphism (\cite{Mie}).
\end{rem}

We say that a non-zero weight
$\lambda\in\Lambda=\Lambda(\lie{u}^\mbb{C},\lie{t})$ is compact, if the
reflection with respect to the hypersurface $(i\eta_\lambda)^\perp\subset
\lie{t}$ belongs to $W$. Otherwise, $\lambda$ is called non-compact. Let us
assume that there exists a good ordering $\Lambda^+$ on the set of non-zero
weights, i.\,e.\ that that each positive non-compact weight is larger than every
compact weight. This implies that the convex cone \begin{equation*}
C_{\max}:=\bigl\{\eta\in\lie{t};\ i\lambda(\eta)\geq0\text{ for all non-compact
$\lambda\in\Lambda^+$}\bigr\}\subset\lie{t}
\end{equation*}
is $W$--invariant.

Let $\Omega:=G_1\exp(iC^0_{\max})G_2$, where $C^0_{\max}$ is the interior of
$C_{\max}$. Then $G_1\times G_2$ acts properly on $\Omega$. Moreover,
the mapping $\cal{R}\colon\cal{C}^\infty(\Omega)^{G_1\times G_2}\to
\cal{C}^\infty(C^0_{\max})^W$,
\begin{equation*}
\Phi\mapsto\varphi:\eta\mapsto\Phi\bigl(\exp(i\eta)\bigr),
\end{equation*}
is an isomorphism (compare~\cite{Fl}). The inverse $\cal{E}:=\cal{R}^{-1}$ is
called the extension operator.

One would expect that the Levi form of an invariant smooth function on $\Omega$
is determined by the direction tangent to the $(G_1\times G_2)$--orbits and by
a direction transversal to the orbit. The following propositions explains how
the Levi form $L(\Phi)(z)$ is influenced by the complex tangent space of
$(G_1\times G_2)\cdot z=M_z$.

\begin{lem}\label{LeviFormofInvFcts}
Let $\Phi\in\cal{C}^\infty(\Omega)^{G_1\times G_2}$ be given. If $v,w\in
H_zM_z\subset T_z\Omega$, then we have
\begin{equation*}
L(\Phi)(z)(v,w)=-d^c\Phi(z)\cal{L}_z(v,w),
\end{equation*}
where $\cal{L}_z$ is the Levi form of $M_z$.
\end{lem}

\begin{proof}
By definition, the Levi form of $\Phi\in\cal{C}^\infty(\Omega)$ at the point
$z\in\Omega$ is the Hermitian form $L(\Phi)(z)$ on $T_z\Omega$ associated to the
$(1,1)$--form $\omega:=-\tfrac{1}{2}dd^c\Phi$. We use the formula
\begin{equation*}
d\omega(V,W)=V\bigl(\omega(W)\bigr)-W\bigl(\omega(V)\bigr)-
\omega\bigl([V,W]\bigr)
\end{equation*}
and extend $v$ to a CR vector field $V$ on $M_z$ to compute as follows:
\begin{align*}
-dd^c\Phi(z)(v,J_zv)&=-v\bigl(d^c\Phi(JV)\bigr)+J_zv\bigl(d^c\Phi(V)\bigr)
+d^c\Phi(z)[V,JV]\\
&=v\bigl(d\Phi(V)\bigr)+J_zv\bigl(d\Phi(JV)\bigr)+d^c\Phi(z)[V,JV]\\
&=v\bigl(V(\Phi)\bigr)+J_zv\bigl(JV(\Phi)\bigr)+d^c\Phi(z)[V,JV].
\end{align*}
Since the vector fields $V$ and $JV$ are tangent to the orbit and since $\Phi$
is constant along the orbit, we obtain
\begin{equation*}
L(\Phi)(z)(v,v)=-\tfrac{1}{2}dd^c\Phi(z)(v,Jv)=\tfrac{1}{2}d\Phi(z)J[V,JV]
_z=-d\Phi(z)J\widehat{\cal{L}}_z(v).
\end{equation*}
Thus the claim follows from the polarization identities.
\end{proof}

\begin{prop}\label{Orthogonality}
Let $z\in\Omega\cap U^\mbb{C}_{sr}$ and let $T_z\Omega=T_zU^\mbb{C}$ be
identified with
\begin{equation}\label{decomp}
\lie{u}^\mbb{C}=\bigoplus_{(\lambda,a)\in\wt{\Lambda}}
\lie{u}^\mbb{C}_{\lambda,a}
\end{equation}
via $(\ell_z)_*$. Let $\varphi\in\cal{C}^\infty(C^0_{\max})^W$ be given and let
$\Phi:=\cal{E}(\varphi)$ be its extension to a smooth $(G_1\times
G_2)$--invariant function on $\Omega$. Then the decomposition~\eqref{decomp} is
orthogonal with respect to the Levi form $L(\Phi)(z)$.
\end{prop}

\begin{proof}
In view of Lemma~\ref{LeviFormofInvFcts} it is enough to show that
$\lie{t}^\mbb{C}$ and $H_zM_z\cong\bigoplus_{(\lambda,a)\not=(0,1)}
\lie{u}^\mbb{C}_{\lambda,a}$ are orthogonal with respect to $L(\Phi)(z)$.
Thus let $v\in\lie{t}$ and $w\in\lie{u}^\mbb{C}_{\lambda,a}$ be given. Since
$Jv$ and $w$ are tangent to $M_z=(G_1\times G_2)\cdot z$, there are elements
$\eta,\xi\in\lie{q}_z\subset\lie{g}_1\oplus\lie{g}_2$ such that
$Jv=\eta_\Omega(z)$ and $w=\xi_\Omega(z)$ hold, where $\eta_\Omega$ and
$\xi_\Omega$ are the corresponding vector fields on $\Omega$. Using the same
arguments as in the proof of Lemma~\ref{LeviFormofInvFcts} together with the
invariance of $\Phi$ we obtain
\begin{equation*}
L(\Phi)(z)(v,w)=d^c\Phi(z)[\eta_\Omega,\xi_\Omega](z)-id^c\Phi(z)
[\eta_\Omega,\xi_\Omega](z).
\end{equation*}
Since $[\eta_\Omega,\xi_\Omega](z)=[\eta,\xi]_\Omega(z)\in H_zM_z$,
the invariance of $\Phi$ implies $d^c\Phi(z)[\eta_\Omega,\xi_\Omega](z)=0$,
which finishes the proof.
\end{proof}

We will apply Proposition~\ref{Orthogonality} in order to establish existence of
a strictly $q$--pseudo-convex exhaustion function on $\Omega$. The following
theorem extends Neeb's result on open complex Ol'shanski{\u\i} semi-groups to
the case $G_1\not=G_2$.

\begin{thm}
The domain $\Omega $ is $q$--complete for
\begin{equation*}
q=\rk(\lie{g}_1\cap\lie{g}_2)+\#\bigl\{(\lambda,1)\in\wt{\Lambda};\
\lambda\in\Lambda^+\bigr\}
+\#\bigl\{(\lambda,-1)\in\wt{\Lambda};\ \lambda\in\Lambda^+\bigr\}
+\#\bigl\{(\lambda,a)\in\wt{\Lambda};\ a\not=\pm1\bigr\}.
\end{equation*}

\end{thm}

\begin{proof}
Let $\varphi\colon C_{\max}^0\to\mbb{R}^{\geq0}$ be smooth, $W$--invariant, and
strictly convex with the property that $\varphi(x_n)\to\infty$ whenever $x_n
\to x\in\partial C_{\max}^0$, and let $\Phi:={\cal{E}}(\varphi)$ be the
corresponding smooth $(G_1\times G_2)$--invariant function in $\Omega$. Let
$z\in\Omega\cap U^\mbb{C}_{sr}$. Due to Proposition~\ref{Orthogonality} we may
compute the Levi form $L(\Phi)(z)$ on each $\lie{u}^\mbb{C}[\lambda,a]$
separately.

We start with the case $a\not=1$. Then our considerations from the
rank one case imply together with Lemma~\ref{LeviFormofInvFcts} that we obtain
for each $\lie{u}^\mbb{C}[\lambda,-1]$ a pair of one positive and one negative
eigenvalue and for each $\lie{u}^\mbb{C}[\lambda,a]$, $a\not=-1$, two pairs of
positive and negative eigenvalues in the Levi form $L(\Phi)(z)$.

Thus let $a=1$. In this case all computations take place in
$(\lie{g}_1\cap\lie{g}_2)^\mbb{C}$ and hence the whole question is reduced to
the case that $\Omega$ is an open Ol'shanski{\u\i} semi-group in $(G_1\cap
G_2)^\mbb{C}$. In~\cite{Ne2} it is proven that in this case the extension of a
strictly convex function is strictly plurisubharmonic. Hence, we see that in our
case the extension $\Phi$ is strictly $q$--pseudo-convex for
\begin{equation*}
q=\dim\lie{t}^\mbb{C}+\#\bigl\{(\lambda,a)\in\wt{\Lambda};\
a=1\bigr\}
+\#\bigl\{(\lambda,-1)\in\wt{\Lambda};\ \lambda\in\Lambda^+\bigr\}
+\#\bigl\{(\lambda,a)\in\wt{\Lambda};\ a\not=\pm1\bigr\}
\end{equation*}
If $z_n\to z\in\partial\Omega$, then $\Phi(z_n)\to\infty$ holds by construction
and hence we conclude that $\Omega$ is $q$--complete for the above $q$.
\end{proof}

\vspace{1cm}

{\setlength{\parindent}{0cm}
{\bf Author's adress:}\\
Fakult\"at f\"ur Mathematik\\
Ruhr-Universit\"at Bochum\\
Universit\"atsstrasse 150\\
44780 Bochum\\
Germany\\

{\bf e-mail:}\\
christian.miebach@ruhr-uni-bochum.de}

\end{document}